\documentclass[leqno,10pt]{article}
\usepackage{amssymb,amsmath}
\usepackage[all]{xy}

\newenvironment{proof}{{\bf Proof. }}{\par}{\bigskip}
\newtheorem{theo}{Theorem}[section]
\newtheorem{defi}[theo]{Definition}

\newtheorem{lem}[theo]{Lemma}
\newtheorem{prop}[theo]{Proposition}

\newtheorem{rem}[theo]{Remark}

\newtheorem{exam}[theo]{Example}

\newcommand{\kgot}{\ensuremath{\mathfrak{k}}}
\newcommand{\qgot}{\ensuremath{\mathfrak{q}}}
\newcommand{\ggot}{\ensuremath{\mathfrak{g}}}
\newcommand{\hgot}{\ensuremath{\mathfrak{h}}}

\newcommand{\Acal}{\ensuremath{\mathcal{A}}}

\newcommand{\Ecal}{\ensuremath{\mathcal{E}}}

\newcommand{\Hcal}{\ensuremath{\mathcal{H}}}

\newcommand{\Jcal}{\ensuremath{\mathcal{J}}}

\newcommand{\Kcal}{\ensuremath{\mathcal{K}}}

\newcommand{\Ncal}{\ensuremath{\mathcal{N}}}

\newcommand{\Scal}{\ensuremath{\mathcal{S}}}
\newcommand{\Ucal}{\ensuremath{\mathcal{U}}}
\newcommand{\Vcal}{\ensuremath{\mathcal{V}}}

\newcommand{\Pbb}{\ensuremath{\mathbb{P}}}

\newcommand{\Cbb}{\ensuremath{\mathbb{C}}}
\newcommand{\Rbb}{\ensuremath{\mathbb{R}}}
\newcommand{\Zbb}{\ensuremath{\mathbb{Z}}}

\newcommand{\A}{\ensuremath{\mathbb{A}}}
\newcommand{\F}{\ensuremath{\mathbf{F}}}

\newcommand{\e}{\operatorname{e}}

\newcommand{\KK}{\ensuremath{{\mathbf K}^0}}

\newcommand{\f}{\ensuremath{\mathcal{C}^{\infty}}}
\newcommand{\fgene}{\ensuremath{\mathcal{C}^{-\infty}}}
\newcommand{\croc}{\ensuremath{\hookrightarrow}}
\newcommand{\T}{\ensuremath{\hbox{\bf T}}}
\newcommand{\End}{\ensuremath{\hbox{\rm End}}}

\newcommand{\str}{\operatorname{Str}}

\newcommand{\tr}{\operatorname{Tr}}
\newcommand{\eul}{\operatorname{Eul}}

\newcommand{\Df}{\operatorname{D}}
\newcommand{\Ag}{\operatorname{\widehat{A}}}

\newcommand{\ch}{\operatorname{Ch}}

\newcommand{\chs}{\operatorname{Ch_{\rm sup}}}

\newcommand{\chg}{\operatorname{Ch_{\rm sup}}}
\newcommand{\chcf}{\operatorname{Ch^1_{\rm c}}}
\newcommand{\chc}{\operatorname{Ch_{\rm c}}}

\newcommand{\supp}{\operatorname{\hbox{\rm \small supp}}}

\renewcommand{\index}{\operatorname{index}}

\def \clif {\mathbf{c}}

\def \bott {{\rm Bott}}

\def \mul {{\rm m}}

\def \Ad {{\rm Ad}}
\def \at {{\rm At}}
\def \ad {{\rm ad}}

\def \spin {{\rm Spin}}

\def \clif {{\bf c}}

\def \Par {{\rm Par}}

\title{Index of transversally elliptic operators}
\author{Paul-Emile Paradan, Mich\`ele Vergne}

\date{April 2008}

\begin{document}

\maketitle

 {\small
 \tableofcontents}

\section{Introduction}

Let $M$ be a compact manifold. The Atiyah-Singer formula for the
index of an {\em elliptic} pseudo-differential operator $P$ on $M$
with {\em elliptic} symbol $\sigma$ on $\T^*M$ involves
integration over the non compact manifold $\T^*M$ of the Chern
character $\ch_c(\sigma)$ of $\sigma$ multiplied by the square of the $\Ag$-genus of $M$:
$$
\index(P)=(2i\pi)^{-\dim M}\int_{\mathbf{T}^*M}\Ag(M)^2\ch_{c}(\sigma).
$$

Here $\sigma$, the principal symbol of $P$,  is a morphism of vector bundles
on $\T^*M$ invertible outside the zero section of $\T^*M$ and the
Chern character $\ch_{c}(\sigma)$ is supported on a small
neighborhood of $M$ embedded in $\T^*M$ as the zero section. It is
important that the representative of the Chern character
$\ch_c(\sigma)$ is compactly supported to perform integration.

Assume that a compact Lie group $K$ (with Lie algebra $\kgot$)
acts on $M$. If the elliptic operator $P$ is  $K$-invariant, then
$\index(P)$ is a smooth function on $K$. The equivariant index of
$P$ can be expressed similarly as the integral of the equivariant
Chern character  of $\sigma$ multiplied by the square of the
equivariant $\Ag$-genus of $M$: for $X\in \kgot$ small enough,
\begin{equation}\label{delocalized}
\index(P)(\e^X)=(2i\pi)^{-\dim M}\int_{\mathbf{T}^*M} \Ag(M)^2(X)\ch_{c}(\sigma)(X).
\end{equation}

Here $\ch_{c}(\sigma)(X)$ is a compactly supported closed equivariant
differential form, that is a differential form on $\T^*M$
depending smoothly of $X\in \kgot$, and closed for the equivariant
differential $D$. The result of the integration determines a
smooth function on a neighborhood of $1$ in $K$ and similar
formulae can be given near any point of $K$. Formula
(\ref{delocalized}) is a ``delocalization" of the
Atiyah-Bott-Segal-Singer formula, in the sense of Bismut \cite{Bismut-localization}.

 The delocalized index formula (\ref{delocalized}) can be adapted to
new cases such as:
\begin{itemize}
\item Index of transversally elliptic operators.

\item $L^2$-index of some elliptic operators  on some non-compact
manifolds (Rossmann formula for discrete series \cite{Rossman78}).
\end{itemize}

Indeed,  in these two contexts, the index exists in the sense of
generalized functions but cannot be always computed in terms of
fixed  point formulae. A ``delocalized" formula will however
continue to have a meaning, as we explain now for transversally
elliptic operators.

  The invariant  operator $P$ with symbol
$\sigma(x,\xi)$ on $\T^*M$ is called transversally elliptic,  if
it is elliptic in the directions transversal to $K$-orbits. In
this case, the operator $P$ has again an index which is a
generalized function on $K$. A very simple example of
transversally elliptic operator is the operator $0$ on $L^2(K)$:
its index is the trace of the action of $K$ in $L^2(K)$, that is
the $\delta$-function on $K$. At the opposite side, $K$-invariant
elliptic operators are of course transversally elliptic, and index
of such operators are smooth functions on $K$ given by Formula
(\ref{delocalized}). Thus a cohomological formula must incorporate
these two extreme cases. Such a cohomological formula was given in
Berline-Vergne   \cite{BV.inventiones.96.1,BV.inventiones.96.2}.
We present here a new point of view, where the  equivariant Chern
character  $\ch_{c}(\sigma)(X)$ entering in Formula
(\ref{delocalized}) is replaced by  a Chern character with
generalized  coefficients, but still {\bf compactly supported}.
Let us briefly explain the construction.

We denote by $\T^*_KM$ the conormal bundle to the $K$-orbits. An
element $(x,\xi)$ is in $\T^*_KM$ if $\xi$ vanishes on tangent
vectors to $K$-orbits. Let $\supp(\sigma)$ be the support of the
symbol $\sigma$ of a transversally elliptic operator $P$. By
definition, the intersection $\supp(\sigma)\cap \T^*_K(M)$ is
compact. By Quillen super-connection construction, the Chern
character $\ch(\sigma)(X)$ is a closed equivariant differential
form supported near the closed set $\supp(\sigma)$. Using the
Liouville one form $\omega$ of $\T^*M$, we construct a closed
equivariant form $\Par(\omega)$ supported  near $\T^*_KM$. Outside
$\T^*_K M$, one has indeed the equation
$1=D(\frac{\omega}{D\omega})$, where the inverse of the form
$D\omega$ is defined by $-i\int_0^{\infty}e^{it D\omega}dt$,
integral which is well defined in the generalized sense, that is
tested against a smooth compactly supported density on $\kgot$.
Thus using a function $\chi$ equal to $1$ on a small neighborhood
of $\T^*_K M$, the closed equivariant form
$$
\Par(\omega)(X)=\chi+d\chi \frac{\omega}{D\omega(X)},\quad X\in\kgot,
$$
is  well defined, supported near $\T^*_KM$, and represents $1$ in
cohomology without support conditions. Remark that
$$\ch_c(\sigma,\omega):=\ch(\sigma)(X)\Par(\omega)(X)$$ is {\bf
compactly supported}. We prove that, for $X\in \kgot$ small
enough, we have

\begin{equation}\label{eq:formule-indice-intro}
\index(P)(\e^X)=(2i\pi)^{-\dim M}\int_{\mathbf{T}^*M}
\Ag(M)^2(X)\ch(\sigma)(X)\,\Par(\omega)(X) .
\end{equation}

 This formula is  thus entirely
similar to the delocalized version of the Atiyah-Bott-Segal-Singer
equivariant index theorem. We have just localized the formula for
the index near $\T^*_K M$ with the help of the form
$\Par(\omega)$, equal to $1$ in cohomology, but supported near
$\T^*_K M$.

 When $P$ is elliptic we can furthermore
localize on the zeros of $VX$ (the vector field on $M$ produced by
the action of $X$) and we obtain the Atiyah-Bott-Segal-Singer
fixed point formulae for the equivariant index of $P$. However the
main difference is that usually we cannot obtain fixed point
formula for the index. For example, the index of a transversally
elliptic operator $P$ where $K$ acts freely is a generalized
function on $K$ supported at the origin. Thus in this case the use
of the form $\Par(\omega)$ is essential. Its role is clearly
explained in the example of the $0$ operator on $S^1$ given at the
end of this introduction.

We  need also to define the  formula for the index at any point of
$s\in K$, in function of integrals over $\T^*M(s)$, where $M(s)$
is the fixed point submanifold of $M$ under the action of $s$. The
compatibility properties (descent method) between the formulae at
different points $s$ are easy to prove, thanks to a localization
formula adapted to this generalized setting.

In the Berline-Vergne 
cohomological formula for the index of  $P$, the Chern character $\ch_c(\sigma)(X)$  in Formula
 (\ref{delocalized}) was replaced by a  Chern character \break
 $\ch_{BV}(\sigma,\omega)(X)$ depending also of the Liouville one form $\omega$.
This Chern character $\ch_{BV}(\sigma,\omega)$ is constructed for
``good symbols" $\sigma$.  It looks like a Gaussian   in the
transverse directions, and is oscillatory in the directions of the
orbits. Our new point of view defines the {\bf compactly supported
} product class $\ch(\sigma)\Par(\omega)$ in a straightforward
way.
 We proved in \cite{pep-vergneIII} that the
classes $\ch_{BV}(\sigma,\omega)$ and $\ch(\sigma)\Par(\omega)$
are equivalent in an appropriate cohomology space, so  that our
new cohomological formula gives the analytic index. However, in
this article, we  choose to prove directly the equality between
the analytic index and the cohomological index, since we want to
show that our formula  in terms of  the  product class
$\ch(\sigma)\Par(\omega)$ is  natural. We follow the same line
than Atiyah-Singer: functoriality with respect to products and
free actions. The compatibility with the free action reduces
basically to the case of the zero operator on $K$, and the
calculation is straightforward. The typical calculation is shown
below. The multiplicativity property is more delicate, but is
based on a general principle on multiplicativity of relative Chern
characters that we proved in a preceding article
\cite{pep-vergneIII}. Thus, following Atiyah-Singer
\cite{Atiyah.74}, we are reduced to the case of $S^1$ acting on a
vector space. The basic examples are then the pushed symbol with
index $-\sum_{n=1}^{\infty}e^{in\theta}$ and the index of the
tangential $\overline \partial$ operators on odd dimensional
spheres. We include at the end a general formula due to
Fitzpatrick \cite{Fitzpatrick07} for contact manifolds.

Let us finally point out that there are many examples of
transversally elliptic operators of great interest. The index of
elliptic  operators on orbifolds  are best understood as indices
of transversally elliptic operators on manifolds where a group $K$
acts with finite stabilizers. The restriction to the maximal
compact subgroup $K$ of a representation of the discrete series of
a real reductive group are indices of transversally elliptic
operators \cite{pep5}. More generally, there is a canonical
transversally elliptic operator on any prequantized  Hamiltonian
manifold with proper moment map (under some mild assumptions)
\cite{pep5},\cite{VergneICM}. For example, when a torus $T$ acts on a Hermitian vector $V$ 
with a proper moment map, then the partition function which computes the $T$-multiplicities on the polynomial algebra $\Scal(V^*)$ is equal 
to the index of a $T$-transversally elliptic operator on $V$. Furthermore, as already noticed  in
Atiyah-Singer, and systematically used in \cite{pep4},
transversally elliptic operators associated to symplectic vector
spaces with proper moment maps and  to cotangent manifolds $\T^*K$
are the local building pieces of any $K$-invariant elliptic
operator.

\begin{exam} Let us check the validity of (\ref{eq:formule-indice-intro}) in the example of the zero operator $0_{S^1}$
acting on the circle group $S^1$. This operator is $S^1$-transversally
elliptic and its index is equal to
$$
\delta_1(\e^{iX})=\sum_{k\in\Zbb}\e^{ikX}, \quad X\in {\rm Lie}(S^1)\simeq \Rbb.
$$

The principal symbol $\sigma$ of $0_{S^1}$ is the zero morphism on
the trivial bundle $S^1\times \Cbb$. Hence $\ch(\sigma)(X)=1$. The
equivariant class $\Ag(S^1)^2(X)$ is also equal to $1$. Thus the
right hand side of (\ref{eq:formule-indice-intro}) becomes
$$
(2i\pi)^{-1}\int_{\mathbf{T}^*S^1}\Par(\omega)(X).
$$
The cotangent bundle $\T^*S^1$ is parametrized by $(\e^{i\theta},\xi)\in S^1\times \Rbb$. The Liouville $1$-form
is $\omega=-\xi d\theta$ : the symplectic form $d\omega=d\theta\wedge d\xi$ gives the orientation of $\T^*S^1$.
Since $VX=-X\frac{\partial}{\partial \theta}$, we have $D\omega(X)= d\theta\wedge d\xi - X\xi$.

Let $g\in \f(\Rbb)$ with compact support and equal to $1$ in a neighborhood of $0$. Then $\chi=g(\xi^2)$ is a
function on $\T^*S^1$ which is supported in a neighborhood of $\T^*_{S^1}S^1={\rm zero\ section}$. We look now at the
equivariant form $\Par(\omega)(X)=\chi+d\chi\wedge(-i\omega)\int_0^{\infty}\e^{it D\omega(X)}dt$. We have
\begin{eqnarray*}
\Par(\omega)(\e^{i\theta},\xi, X)
&=& g(\xi^2)+  g'(\xi^2)2\xi d \xi\wedge(i\xi d\theta) \int_0^{\infty}\e^{it (d\theta\wedge d\xi - X\xi)}dt\\
&=& g(\xi^2)-   i d\theta \wedge d(g(\xi^2))\left(\int_0^{\infty}\e^{-itX\xi}\xi dt \right) .
\end{eqnarray*}
If we make the change of variable $t\xi\to t$ in the integral $\int_0^{\infty}\e^{-itX\xi}\xi dt$ we get
$$
\Par(\omega)(\e^{i\theta},\xi, X)= \begin{cases}
  g(\xi^2)  -   i d\theta\wedge d(g(\xi^2))\left(\int_0^{\infty}\e^{-itX}dt \right) , & \text{if}\quad \xi \geq 0;\\
   g(\xi^2) +  i d\theta\wedge d(g(\xi^2))\left(\int^0_{-\infty}\e^{-itX}dt \right) , & \text{if}\quad \xi\leq  0.
\end{cases}
$$
Finally, since $-\int_{\xi\geq 0}d(g(\xi^2))=\int_{\xi\leq 0}d(g(\xi^2))= 1$, we have
$$
(2i\pi)^{-1}\int_{\mathbf{T}^*S^1}\Par(\omega)(X)= \int_{-\infty}^{\infty}\e^{-itX}dt.
$$
The generalized function $\delta_0(X)=\int_{-\infty}^{\infty}\e^{-itX}dt$ satisfies
$$
\int_{{\rm Lie}(S^1)}\delta_0(X)\varphi(X)dX= {\rm vol}(S^1,dX) \varphi(0)
$$
for any function $\varphi\in \f({\rm Lie}(S^1))$ with compact support. Here ${\rm vol}(S^1,dX)=\int_{0}^{2\pi}dX$ is also
the volume of $S^1$ with the Haar measure compatible with $dX$.

Finally, we see that (\ref{eq:formule-indice-intro}) corresponds to the following equality of generalized functions
$$
\delta_1(\e^{iX})=\delta_0(X),
$$
which holds for $X\in {\rm Lie}(S^1)$ small enough.

\end{exam}

\section{The analytic index}

\subsection{Generalized functions}

Let $K$ be a compact Lie group. We denote by  $\hat K$
the set of unitary irreducible representations of $K$. If $\tau\in
\hat K$, we denote by $V_\tau$ the representation space of $\tau$ and by $k\mapsto\tr(k,\tau)$
its character. Let $\tau^*$ be the dual representation of $K$ in $V_\tau^*$.

We denote by $C^{-\infty}(K)$ the space of generalized
functions on $K$ and by $C^{-\infty}(K)^K$ the space of central invariant
generalized functions. The space $C^{\infty}(K)$ of
smooth functions on $K$ is naturally a subspace of
$C^{-\infty}(K)$. We will often use the notation  $\Theta(k)$ to
denote a generalized function $\Theta$ on $K$, although (in
general) the value of $\Theta$ on a particular point $k$ of $K$
does not have a meaning. By definition, $\Theta$ is a linear form
on the space of smooth densities on $K$. If $dk$ is a Haar measure
on $K$ and $\Phi\in C^{\infty}(K)$, we denote by $\int_K
\Theta(k)\Phi(k) dk$ the value of $\Theta$ on the density $\Phi
dk$.

Any invariant generalized function on $K$ is expressed as
$\Theta(k)=$ \break $\sum_{\tau\in \hat K}n_\tau \tr(k,\tau)$ where the Fourier
coefficients $n_\tau$ have at most of polynomial growth.

\subsection{Symbols and pseudo-differential operators}

Let $M$ be a compact manifold with a smooth action of a compact
Lie group $K$. We consider the closed subset $\T_K^*M$ of the
cotangent bundle $\T^*M$, union of the spaces $(\T^*_KM)_x,x\in
M$, where  $(\T^*_KM)_x\subset \T^*_x M$ is the orthogonal of the
tangent space at $x$ to the orbit $K\cdot x$. Let $\Ecal^{\pm}$ be two
$K$-equivariant complex vector bundles over $M$. We denote by
$\Gamma(M,\Ecal^\pm)$ the space of smooth sections of $\Ecal^\pm$.
Let $P:\Gamma(M,\Ecal^+)\to\Gamma(M,\Ecal^-)$ be a $K$-invariant
pseudo-differential operator of order $m$. Let $p:\T^*M\to M$ be
the natural projection. The principal symbol $\sigma(P)$ of $P$ is
a bundle map $p^*\Ecal^+\to p^*\Ecal^-$ which is homogeneous of
degree $m$, defined over $\T^*M\setminus M$.

The operator $P$ is elliptic if its principal symbol $\sigma(P)(x,\xi)$ is invertible for all
$(x,\xi)\in \T^* M$ such that $\xi\neq 0$. The operator $P$ is said to be $K$-transversally
elliptic if its principal symbol $\sigma(P)(x,\xi)$ is invertible for all $(x,\xi)\in \T^*_K M$
such that $\xi\neq 0$.

Using  a $K$- invariant function $\chi$ on $\T^*M$ identically
equal to $1$ in a neighborhood of $M$ and compactly supported,
then $\sigma_P(x,\xi):=(1-\chi(x,\xi)) \sigma(P)(x,\xi)$ is a morphism  from
$p^*\Ecal^+$ to $p^*\Ecal^-$ defined on the whole space $\T^*M$
and which is {\em almost homogeneous}: $\sigma_P(x,t\xi)= t^m\sigma_P(x,\xi)$
for $t>1$ and $\xi$ large enough. We consider the support of the morphism $\sigma_P$,
$$
\supp(\sigma_P):=\{(x,\xi)\in\T^*M\ \vert\  \sigma_P(x,\xi) \ \hbox{is not invertible}\}
$$
which is a closed $K$-invariant subset of $\T^*M$.

When $P$ is elliptic, then $\supp(\sigma_P)$ is compact, and the morphism $\sigma_P$ gives rise
to a $\KK$-theory class $[\sigma_P]\in \KK_K(\T^*M)$ which does not depend of the choice of
$\chi$. Similarly, when $P$ is $K$-transversally elliptic, then $\supp(\sigma_P)\cap \T^*_K M$
is compact and the morphism $\sigma_P$ gives rise to a $\KK$-theory class
$[\sigma_P\vert_{\mathbf{T}^*_K M}]\in \KK_K(\T^*_K M)$ which does not depend of
the choice of $\chi$.

Recall the definition of the $K$-equivariant index of a pseudo-differential operator $P$
which is $K$-transversally elliptic. Let us choose a $K$-invariant metric on $M$ and
$K$-invariant Hermitian structures on $\Ecal^{\pm}$. Then the adjoint $P^*$ of $P$ is also
a $K$-transversally elliptic pseudo-differential operator.

If $P$ is elliptic, its kernel $\ker P:=\{s\in \Gamma(M,\Ecal^+)|\; Ps=0 \}$ is finite dimensional,  and
the $K$-equivariant index of $P$ is the central function $\index^K(P)(k)=\tr(k,\ker P)-\tr(k, \ker P^*)$.

If $P$ is $K$-transversally elliptic, its kernel $\ker P$ is not
finite dimensional, but it has finite multiplicities: the vector
space $\ker P$ is an admissible $K$-representation. Let us explain
this notion.  For any irreducible representation $\tau\in\hat{K}$,
the multiplicity $m_\tau(P):= \hom_K(V_\tau, \ker P)$ is finite,
and $\tau\mapsto m_\tau(P)$ has at most a polynomial growth. We
define then a central invariant generalized function on $K$ by
setting
$$
\tr(k,\ker P):= \sum_{\mu\in\hat{K}}m_\tau(P)\tr(k,\tau).
$$

\begin{defi}
The $K$-equivariant index of a $K$-transversally elliptic pseudo-differential operator $P$
is the generalized function
$$\index^K(P)(k)=\tr(k,\ker P)-\tr(k, \ker P^*).$$
\end{defi}

We recall
\begin{theo}[Atiyah-Singer]
$\bullet$ The $K$-equivariant index of a $K$-invariant elliptic pseudo-differential operator
$P$ depends only of $[\sigma_P]\in\KK_K(\T^*M)$.

$\bullet$ The $K$-equivariant index of a $K$-transversally elliptic pseudo-differential operator
$P$ depends only of $[\sigma_P\vert_{\mathbf{T}^*_K M}]\in\KK_K(\T_K^*M)$.

$\bullet$ Each element in $\KK_K(\T_K^*M)$ is represented by the class
$[\sigma_P\vert_{\mathbf{T}^*_K M}]$ of a $K$-transversally elliptic
pseudo differential operator $P$ of order $m$. Similarly, each element in $\KK_K(\T^*M)$
is represented  by the class
$[\sigma_P]$ of a $K$-invariant elliptic pseudo differential operator $P$ of order $m$.
\end{theo}

Thus we can define
\begin{equation}\label{eq:def-a-index}
\index_a^{K,M}: \KK_K(\T^*_K M)\to C^{-\infty}(K)^{K}
\end{equation}
by setting,  for $P$ a $K$-transversally elliptic pseudo-differential operator of
order $m$, $\index_a^{K,M}([\sigma_P\vert_{\mathbf{T}^*_K M}])=\index^K(P)$.  Similarly,
we can define in the elliptic setting
\begin{equation}\label{eq:def-a-index-elliptic}
\index_a^{K,M}: \KK_K(\T^* M)\to \f(K)^{K}.
\end{equation}
Note that we have a natural restriction map $\KK_K(\T^* M)\to\KK_K(\T^*_K M)$ which make
the following diagram
\begin{equation}\label{indice.generalise}
\xymatrix{
\KK_{K}(\T^* M)\ar[r]\ar[d]_{\index_a^{K,M}} &
\KK_{K}(\T_{K}^*M)\ar[d]^{\index_a^{K,M}}\\
\f(K)^{K}\ar[r] &  \fgene(K)^{K}
   }
\end{equation}
commutative.

Let $R(K)$ be the representation ring of $K$. Using the trace, we will consider $R(K)$ as a sub-ring
of  $\f(K)^{K}$. The map (\ref{eq:def-a-index}) and (\ref{eq:def-a-index-elliptic}) are homomorphism
of $R(K)$-modules and will be called the analytic indices.

\begin{rem}
In order to simplify the notations we will make no distinction
between an element $V\in R(K)$, and its trace function $k\to
\tr(k,V)$ which belongs to $\f(K)^{K}$. For example the constant
function $1$ on $K$ is identified with the trivial representation
of $K$.
\end{rem}

\bigskip

Let $H$ be a compact Lie group acting on $M$ and commuting with
the action of $K$. Then the space $\T^*_KM$ is provided with an
action of $K\times H$. If $[\sigma]\in \KK_{K\times H}(\T^*_KM)$,
we can associate to $[\sigma]$ a virtual trace class
representation of $K\times H$. Indeed, we can choose as
representative of $[\sigma]$ the symbol of a $H$-invariant and
$K$-transversally elliptic operator $P$. Then $\ker P-\ker P^*$
is a trace class representation of $K\times H$. Thus we can define
a $R(K\times H)$-homomorphism:
$$\index^{K,H,M}_a:\KK_{K\times H}(\T^*_KM)\to C^{-\infty}(K\times
H)^{K\times H}.$$

Obviously $\T^*_{K\times H}(M)$ is contained in $\T^*_K(M)$ so we have a restriction
morphism $r: \KK_{K\times H}(\T^*_KM)\to \KK_{K\times H}(\T^*_{K\times H} M)$. We see
that
$$\index^{K,H,M}_a=\index^{K\times H,M}_a\circ\  r.
$$

However, it is easy to see that for a $H$-invariant and
$K$-transversally elliptic symbol $[\sigma]$,
$\index^{K,H,M}_a([\sigma])$ belongs to
$C^{\infty}(H,C^{-\infty}(K))$. In particular we can restrict
$\index^{K,H,M}_a([\sigma])$ to $H'\times
K$, for a subgroup $H'$ of $H$. We can also multiply
$\index^{K,H,M}_a([\sigma])(g,h)$ by generalized functions
$\Psi(h)$ on $H$.

\subsection{Functoriality property of the analytic index}

We have defined for any compact $K\times H$-manifold $M$ a $R(K\times H)$-morphism
$$
\index^{K,H,M}_a:\KK_{K\times H}(\T^*_KM)\to \f(H,C^{-\infty}(K))^{K\times H}.
$$

Let us recall some basic properties of the analytic index map:
\begin{itemize}
\item $\mathbf{[N1]}$ If $M=\{{\rm point}\}$, then $\index^{K,M}_a$ is the trace map
$R(K)\croc \f(K)^{K}$.
\item $\mathbf{[Diff]}$ Compatibility with diffeomorphisms: if $f:M_1\to M_2$ is a
$K\times H$-diffeomorphism then $\index^{K,H,M_1}_a\circ f^*$ is equal to $\index^{K,H,M_2}_a$.
\item $\mathbf{[Morph]}$ If $\phi: H'\to H$ is a Lie group morphism, we have
$\phi^*\circ \ \index^{K,H,M}_a$ $=\index^{K,H',M}_a$.
\end{itemize}

\subsubsection{Excision}
Let $U$ be a non-compact $K$-manifold. Lemma 3.6 of \cite{Atiyah.74} tell us that, for any open
$K$-embedding $j:U\croc M$ into a compact manifold, we have a pushforward map
$j_{*}:\KK_{K}(\T_{K}^* U)\to \KK_{K}(\T_{K}^* M)$.

Let us rephrase Theorem 3.7 of \cite{Atiyah.74}.
\begin{theo}[Excision property]
The composition
$$
\KK_{K}(\T_{K}^*U)\stackrel{j_{*}}{\longrightarrow} \KK_{K}(\T_{K}^*M)
\stackrel{\index^{K,M}_a}{\longrightarrow} \fgene(K)^{K}
$$
is independent of the choice of $j:U\croc M$: we denote this map $\index^{K,U}_a$.
\end{theo}

Note that a {\em relatively compact} $K$-invariant open subset $U$ of a
$K$-manifold admits an open $K$-embedding $j:U\croc M$ into a compact
$K$-manifold. So the index map $\index^{K,U}_a$ is defined in this case.
An important example is when $U\to N$ is a $K$-equivariant vector bundle
over a compact manifold $N$ : we can imbed $U$ as an open subset
of the real projective bundle $\Pbb(U\oplus \Rbb)$.

\subsubsection{Exterior product}\label{subsec:exterior-product}

Let us recall the multiplicative property of the analytic index for the product
of manifolds.  Consider a compact Lie group $K_2$ acting on two
manifolds $M_1$ and $M_2$, and assume that another compact
Lie group $K_1$
acts on $M_1$ commuting with the action of $K_2$.

The external product of complexes on $\T^*M_1$ and $\T^*M_2$ induces
a multiplication (see \cite{Atiyah.74}):
$$
\odot_{\rm ext}: \KK_{K_1\times K_2}(\T^*_{K_1} M_1)\times \KK_{K_2}(\T^*_{K_2} M_2)
\longrightarrow \KK_{K_1\times K_2}(\T^*_{K_1\times K_2} (M_1\times M_2)).
$$

Let us recall the definition of this external product. For $k=1,2$, we consider equivariant
morphisms\footnote{In order to simplify the notation, we do not make the distinctions between
vector bundles on $\T^*M$ and on $M$.}
 $\sigma_k:\Ecal^+_k\to \Ecal_k^-$  on $\T^*M_k$. We consider the equivariant morphism on
$\T^*(M_1\times M_2)$
$$
\sigma_1\odot_{\rm ext}\sigma_2: \Ecal_1^{+}\otimes \Ecal_2^{+}\oplus \Ecal_1^{-}\otimes \Ecal_2^{-}
\longrightarrow \Ecal_1^{-}\otimes \Ecal_2^{+} \oplus \Ecal_1^{+}\otimes \Ecal_2^{-}
$$
defined by
\begin{equation}\label{eq:produit.externe}
\sigma_{1}\odot_{\rm ext} \sigma_{2}=
\left(
\begin{array}{cc}
  \sigma_{1}\otimes {\rm Id} & -{\rm Id}\otimes \sigma_{2}^{*}\\
{\rm Id} \otimes \sigma_{2} & \sigma_{1}^{*}\otimes {\rm Id}
\end{array}
\right)\ .
\end{equation}

We see that the set $\supp(\sigma_{1}\odot \sigma_{2})\subset
\T^*M_1\times\T^*M_2$ is equal to
$\supp(\sigma_{1})\times \supp(\sigma_{2})$.

We suppose now that the morphisms $\sigma_k$ are respectively
$K_k$-transversally elliptic. Since $\T^*_{K_1\times K_2}
(M_1\times M_2)\neq \T^*_{K_1} M_1\times\T^*_{K_2} M_2$, the
morphism $\sigma_{1}\odot_{\rm ext} \sigma_{2}$ is not necessarily
$K_1\times K_2$-transversally elliptic. Nevertheless, if
$\sigma_2$ is taken {\em almost homogeneous} of order $m=0$, then
the morphism $\sigma_{1}\odot_{\rm ext} \sigma_{2}$ is $K_1\times
K_2$-transversally elliptic (see Lemma 4.9 in
\cite{pep-vergneIII}). So the exterior product $a_1\odot_{\rm
ext}a_2$ is the $\KK$-theory class defined by $\sigma_1\odot_{\rm
ext}\sigma_2$, where $a_k=[\sigma_k]$ and $\sigma_2$ is taken
almost homogeneous of order $m=0$.

\begin{theo}[Multiplicative property] \label{theo:multiplicative-property}
For any $[\sigma_1]\in \KK_{K_1\times K_2}(\T^*_{K_1} M_1)$ and
any $[\sigma_2]\in\KK_{K_2}(\T^*_{K_2} M_2)$ we have
$$
\index^{K_1\times K_2,M_1\times M_2}_a([\sigma_1]\odot_{\rm ext}[\sigma_2])
=\index^{K_1,K_2,M_1}_a([\sigma_1])\index^{K_2,M_2}_a([\sigma_2]).
$$
\end{theo}

\subsubsection{Free action}\label{subsec:a-index-free}

Let $K$ and $G$ be two compact Lie groups. Let $P$ be a compact manifold provided with an action of
$K\times G$. We assume that the action of $K$ is free. Then the manifold $M:=P/K$ is provided with an
action of $G$ and the quotient map $q:P\to M$ is $G$-equivariant. Note that we have the natural
identification of
$\T^*_K P$ with $q^*\T^* M$, hence $(\T^*_K P)/K\simeq \T^*M$ and more generally
$$
(\T^*_{K\times G} P)/K\simeq \T_G^*M.
$$
This isomorphism induces an isomorphism
$$
Q^*: \KK_{G}(\T^*_{G} M)\to \KK_{K\times G}(\T^*_{K\times G} P).
$$
Let $\Ecal^\pm$ be two $G$-equivariant complex vector bundles on $M$ and $\sigma:p^*\Ecal^+\to p^*\Ecal^-$
be a $G$-transversally elliptic symbol. For any irreducible
representation $(\tau,V_\tau)$ of $K$, we form the $G$-equivariant complex vector bundle
$\Vcal_\tau:= P\times_K V_\tau$ on $M$. We consider the morphism
$$
\sigma_\tau:=\sigma\otimes {\rm Id}_{V_\tau}:p^*(\Ecal^+\otimes \Vcal_\tau)\to p^*(\Ecal^-\otimes \Vcal_\tau)
$$
which is $G$-transversally elliptic.

\begin{theo}[Free action property] We have the following equality in \break
$\fgene(K\times G)^{K\times G}$: for $(k,g)\in K\times G$
$$
\index^{K\times G,P}_a(Q^*[\sigma])(k,g)=\sum_{\tau\in \hat{K}} \tr(k,\tau)
\index^{G,M}_a([\sigma_{\tau^*}])(g).
$$
\end{theo}

\subsection{Basic examples}

\subsubsection{Bott symbols}

Let $W$ be a Hermitian vector space. For any $v\in W$, we consider on the
$\Zbb_2$-graded vector space
$\wedge W$ the following odd operators:
the exterior multiplication $\mul(v)$ and the contraction $\iota(v)$. The contraction
$\iota(v)$ is an odd derivation of $\wedge W$ such that
$\iota(v)w=(w,v)$ for $w\in\wedge^1 W=W$.

The Clifford action of $W$ on $\wedge W$ is defined by the formula
\begin{equation}
\clif(v)=\mul(v) - \iota(v)
\end{equation}
Then $\clif(v)$ is an odd operator on $\wedge W$ such that $\clif(v)^2=-\|v\|^2 Id$. So $\clif(v)$ is
invertible when $v\neq 0$.

Consider the trivial vector bundles $\Ecal^\pm := W\times \wedge^\pm W$  over $W$ with
fiber $\wedge^\pm W$.
The Bott morphism $\bott(W): \Ecal^+\to\Ecal^-$ is defined by
\begin{equation}\label{eq:bott-symbol}
\bott(W)(v,w)=(v,\clif(v)w).
\end{equation}

Consider now an Euclidean vector space $V$. Then its complexification $V_\Cbb$ is an Hermitian
vector space. The cotangent bundle $\T^*V$ is identified with $V_\Cbb$: we associate
to the covector $\xi\in \T^*_v V$ the element $v+i\hat{\xi}\in V_\Cbb$, where
$\xi\in V^*\to\hat{\xi}\in V$ is the identification given by the Euclidean structure.

Then $\bott(V_\Cbb)$ defines an elliptic symbol on $V$ which is equivariant relatively to the action of
the orthogonal group $O(V)$.

\begin{prop}
We have $\mathbf{[N2]}$: $\qquad \index^{O(V),V}_a(\bott(V_\Cbb))=1$.
\end{prop}

\begin{rem}\label{rem:reduction-bott-symbol}
If $V$ and $W$ are two euclidean vector spaces we see that the symbol $\bott((V\times W)_\Cbb)$ is
equal to the product $\bott(V_\Cbb) \odot\bott(W_\Cbb)$. Then for $(g,h)\in O(V)\times O(W)$,
the multiplicative property tells us that
$$\index^{O(V\times W),V\times W}_a(\bott((V\times W)_\Cbb))(g,h)$$ is equal
to the product $\index^{O(V),V}_a(\bott(V_\Cbb))(g)\index^{O(W),W}_a(\bott(W_\Cbb))(h)$.

For any $g\in O(V)$, the vector space $V$ decomposes as an orthogonal sum $\oplus_i V_i$ of
$g$-stables subspaces, where either  $\dim V_i=1$ and $g$ acts on $V_i$ as $\pm 1$, or
$\dim V_i=2$ and $g$ acts on $V_i$ as  a rotation.

Hence $\mathbf{[N2]}$ is  satisfied for any euclidean vector space if one check it for the cases:

$\bullet$ $V=\Rbb$ with the action of the group $O(V)=\Zbb_2$,

$\bullet$ $V=\Rbb^2$ with the action of the group $SO(V)=S^1$.

\end{rem}

\subsubsection{Atiyah symbol}\label{subsec:Atiyah-analytic}

In the following example, we denote,
 for any integer $k$, by
$\Cbb_{[k]}$ the vector space $\Cbb$ with the action of the circle
group $S^1$ given by : $t\cdot z=t^k z$.

The Atiyah symbol is the $S^1$-equivariant morphism on
$N=\T^*\Cbb_{[1]}\simeq \Cbb_{[1]}\times\Cbb_{[1]}$
\begin{eqnarray*}
\sigma_\at:N\times \Cbb_{[0]}&\longrightarrow& N\times \Cbb_{[1]}\\
\big(\xi,v\big)&\longmapsto& \big(\xi,\sigma_\at(\xi)v\big)
\end{eqnarray*}
defined by $\sigma_\at(\xi)= \xi_2+i\xi_1$ for $\xi=(\xi_1,\xi_2)\in\T^*\Cbb_{[1]}$.

The symbol $\sigma_\at$ is not elliptic since $\supp(\sigma_\at)=\{\xi_1=i\xi_2\}\subset \Cbb^2$
is not compact. But $\T^*_{S^1}\Cbb_{[1]}=\{(\xi_1,\xi_2)\ \vert\ {\rm Im}\,(\xi_1\overline{\xi_2})=0\}$
and $\supp(\sigma_\at)\cap \T^*_{S^1}\Cbb_{[1]}=\{(0,0)\}$ : the symbol $\sigma_\at$ is
$S^1$-transversally elliptic.
\begin{prop}
We have $\mathbf{[N3]}$: $\qquad\index^{S^1,\Cbb}_a(\sigma_\at)(t)=-\sum_{k=1}^\infty t^k$.
\end{prop}

\subsection{Unicity of the index}

Suppose that for any compact Lie groups $K$ and $H$, and any compact $K\times H$-manifold $M$,
we have a map of $R(K\times H)$-modules:
$$
\mathbb{I}^{K,H,M}:\KK_{K\times H}(\T^*_KM)\to \f(H,C^{-\infty}(K))^{K\times H}.
$$

\begin{theo}
Suppose that the maps $\mathbb{I}^{-}$ satisfy

$\bullet$ the normalization conditions $\mathbf{[N1]},\mathbf{[N2]}$ and $\mathbf{[N3]}$,

$\bullet$ the functorial properties $\mathbf{Diff}$ and $\mathbf{Morph}$,

$\bullet$ the ``excision property'', the ``multiplicative property'' and the ``free action property''.

Then $\mathbb{I}^{-}$ coincides with the analytic index map $\index_a^{-}$.
\end{theo}

\section{The cohomological index}

Let $N$ be a manifold, and let $\Acal(N)$ be the algebra of
differential forms on $N$. We denote by $\Acal_c(N)$ the subalgebra
of compactly supported differential forms. We will
consider on $\Acal(N)$ and $\Acal_c(N)$ the $\Zbb_2$-grading in even
or odd differential forms.

Let $K$ be a  compact Lie group with Lie algebra  $\kgot$.
We suppose that the manifold $N$ is provided with an action of $K$.
We denote $X\mapsto VX$ the corresponding morphism from $\kgot$ into
the Lie algebra of vectors fields on $N$: for $n\in N$,
$V_nX:=\frac{d}{d\epsilon} \exp(-\epsilon X)\cdot n|_{\epsilon=0}$.

Let $\Acal^{\infty}(\kgot, N)$ be the $\Zbb_2$-graded algebra of
equivariant smooth functions $\alpha: \kgot\to\Acal(N)$. Its
$\Zbb_2$-grading is the grading induced by the exterior degree. Let
$D= d-\iota(VX)$ be the equivariant differential:
$(D\alpha)(X)=d(\alpha(X))-\iota(VX)\alpha(X)$. Here the operator
$\iota(VX)$ is the contraction of a differential form by the vector field $VX$.
Let $\Hcal^{\infty}(\kgot,N):= \mathrm{Ker} D/ \mathrm{Im} D$ be the
equivariant cohomology algebra with $C^{\infty}$-coefficients. It is a module
over the algebra $\f(\kgot)^K$ of $K$-invariant $C^{\infty}$-functions on $\kgot$.

The sub-algebra $\Acal^{\infty}_c(\kgot, N)\subset\Acal^{\infty}(\kgot, N)$ of
equivariant differential forms with compact support is defined as follows :
$\alpha\in \Acal^{\infty}_c(\kgot, N)$ if there exists a compact subset $\Kcal_\alpha\subset N$
such that the differential form $\alpha(X)\in \Acal(N)$ is supported on $\Kcal_\alpha$ for
 any $X\in\kgot$. We denote $\Hcal_c^{\infty}(\kgot, N)$ the corresponding algebra of
cohomology: it is a  $\Zbb_2$-graded algebra.

\medskip

Let $\Acal^{-\infty}(\kgot, N)$ be the space of \emph{generalized}
equivariant differential forms. An element $ \alpha\in
\Acal^{-\infty}(\kgot, N)$ is, by definition, a
$\fgene$-equivariant map $\alpha:\kgot\to \Acal(N)$. The value
taken by $\alpha$ on a smooth compactly supported density $Q(X)dX$
on $\kgot$ is denoted by $ \int_{\kgot}\alpha(X)Q(X)dX\in \,
\Acal(N)$. We have $\Acal^{\infty}(\kgot,
N)\subset\Acal^{-\infty}(\kgot, N)$ and we can extend the
differential $D$ to $\Acal^{-\infty}(\kgot, N)$
\cite{Kumar-Vergne}. We denote by $\Hcal^{-\infty}(\kgot,N)$ the
corresponding cohomology space.
Note that $\Acal^{-\infty}(\kgot, N)$ is a module over
$\Acal^{\infty}(\kgot, N)$ under the wedge product, hence the
cohomology space $\Hcal^{-\infty}(\kgot,N)$ is a module over
$\Hcal^{\infty}(\kgot, N)$.

The sub-space $\Acal^{-\infty}_c(\kgot, N)\subset\Acal^{-\infty}(\kgot, N)$ of
generalized equivariant differential forms with compact support is defined as follows :
$\alpha\in \Acal^{-\infty}_c(\kgot, N)$ if there exits a compact subset $\Kcal_\alpha\subset N$
such that the differential form $\int_{\kgot}\alpha(X)Q(X)dX\in \,
\Acal(N)$ is supported on $\Kcal_\alpha$ for
 any  compactly supported density $Q(X)dX$.
We denote $\Hcal_c^{-\infty}(\kgot, N)$ the corresponding space of
cohomology. The $\Zbb_2$-grading on $\Acal(N)$
induces a $\Zbb_2$-grading on the cohomology spaces
$\Hcal^{-\infty}(\kgot,N)$ and $\Hcal^{-\infty}_c(\kgot,N)$.

If $\Ucal$ is a $K$-invariant open subset of $\kgot$, ones defines also
$\Hcal^{-\infty}(\Ucal,N)$ and $\Hcal^{-\infty}_c(\Ucal,N)$. If $N$ is equipped with a $K$-invariant
orientation, the integration over $N$ defines a morphism
$$
\int_N:\Hcal^{-\infty}_c(\Ucal,N)\longrightarrow\fgene(\Ucal)^K.
$$

\subsection{Restrictions of generalized functions}\label{sec:restriction}

Let $K$ be a compact Lie group with Lie algebra $\kgot$. For any $s\in K$ (resp. $S\in \kgot$),
we denote $K(s)$ (resp. $K(S)$) the stabilizer subgroup: the corresponding Lie algebra is denoted
$\kgot(s)$ (resp. $\kgot(S)$).

For any $s\in K$, we consider a (small) open $K(s)$-invariant
neighborhood $\Ucal_s$ of $0$ in $\kgot(s)$ such that the map  $[k,Y]\mapsto kse^Yk^{-1}$ is an open
embedding of
$K\times_{K(s)}\Ucal_s$ on an open neighborhood of the conjugacy
class $K\cdot s:=\{ksk^{-1},\ k\in K\}\simeq K/K(s)$.

Similarly, for any $S\in\kgot$, we consider a (small) open $K(S)$-invariant
neighborhood $\Ucal_S$ of $0$ in $\kgot(S)$ such that the map  $[k,Y]\mapsto \Ad(k)(S+Y)$ is a open
embedding of
$K\times_{K(S)}\Ucal_S$ on an open neighborhood of the adjoint orbit $K\cdot S\simeq K/K(S)$.

Note that the map $Y\mapsto [1,Y]$ realizes $\Ucal_s$ (resp. $\Ucal_S$) as
a $K(s)$-invariant sub-manifold of $K\times_{K(s)}\Ucal_s$ (resp. $K\times_{K(S)}\Ucal_S$).

Let $\Theta$ be a generalized function on $K$ invariant by conjugation. For any $s\in K$, $\Theta$
defines a $K$-invariant generalized function on $K\times_{K(s)}\Ucal_s\croc K$ which admits a
restriction to the submanifold
$\Ucal_s$ that we denote
$$
\Theta\vert_s\in\fgene(\Ucal_s)^{K(s)}.
$$
If $\Theta$ is smooth we have $\Theta\vert_s(Y)=\Theta(s\e^Y)$.

Similarly, let $\theta$ be a $K$-invariant generalized function on $\kgot$. For any $S\in\kgot$,
$\theta$ defines a $K$-invariant generalized function on $K\times_{K(S)}\Ucal_S\croc\kgot$ which
admits a restriction to the submanifold
$\Ucal_S$ that we denote
$$
\theta\vert_S\in\fgene(\Ucal_S)^{K(S)}.
$$
If $\theta$ is smooth we have $\theta\vert_S(Y)=\Theta(S+Y)$.

\medskip

We have $K(se^S)=K(s)\cap K(S)$ for any $S\in \Ucal_s$. Let $\Theta\vert_s\in\fgene(\Ucal_s)^{K(s)}$
be the restriction of a generalized function $\Theta\in\fgene(K)^{K}$. For any $S\in\Ucal_s$, the
generalized function $\Theta\vert_s$ admits a restriction $(\Theta\vert_s)\vert_S$ which is a
$K(s\e^S)$-invariant generalized function
defined in a neighborhood of $0$ in $\kgot(s)\cap \kgot(S)=\kgot(s\e^S)$.

\begin{lem}
Let $\Theta\in\fgene(K)^K$.

$\bullet$ For $s\in K$, and $S\in \Ucal_s$, we have the following equality
of generalized functions defined in a neighborhood of $0$ in $\kgot(s\e^S)$
\begin{equation}\label{eq:compatible}
(\Theta\vert_s)\vert_S=\Theta\vert_{s\e^S}.
\end{equation}

$\bullet$ Let $s, k\in K$. We have the following equality
of generalized functions defined in a neighborhood of $0$ in $\kgot(s)$
\begin{equation}\label{eq:invariante}
\Theta\vert_s=\Theta\vert_{ksk^{-1}}\circ\Ad(k).
\end{equation}
\end{lem}

When $\Theta\in \f(K)^K$ is smooth, condition (\ref{eq:compatible}) is
easy to check: for $Y\in\kgot(s\e^S)$, we have
$$
(\Theta\vert_s)\vert_S(Y)=\Theta\vert_s(S+Y)=\Theta(s\e^{S+Y})=\Theta(s\e^S\e^Y)=\Theta\vert_{s\e^S}(Y).
$$
We have the following

\begin{theo}\label{theo:germe-function}
Let $K$ be a compact Lie group. Consider a family of
generalized function $\theta_s\in
\fgene(\Ucal_s)^{K(s)}$. We assume that the following
conditions are verified.

$\bullet$ {\bf Invariance}: for any $k$ and $s\in  K$, we have the following equality
of generalized functions defined in a neighborhood of $0$ in $\kgot(s)$
$$\theta_s=\theta_{ksk^{-1}}\circ\Ad(k).$$

$\bullet$ {\bf Compatibility}: for every  $s\in  K$ and $S\in U_s$,  we have the following equality
of generalized functions defined in a neighborhood of $0$ in $\kgot(s\e^S)$
$$\theta_s\vert_S=\theta_{s\e^S}.$$

Then there exists a unique  generalized function $\Theta\in
C^{-\infty}(K)^K$ such that, for any $s\in K$,  the equality
$\Theta\vert_s=\theta_s$ holds in $\fgene(\Ucal_s)^{K(s)}$.
\end{theo}

\subsection{Integration of bouquet of equivariant forms}\label{subsec:integration-bouquet}

Let $K$ be a compact Lie group acting on a compact manifold $M$. We are interested in the
central functions on $K$ that can be defined by integrating equivariant forms on $\T^*M$.

Let $\omega$ be the Liouville $1$-form on $\T^*M$.  For any $s\in K$, we denote $M(s)$ the
fixed points set $\{x\in M\ \vert sx=x\}$. As $K$ is compact, $M(s)$ is a submanifold of $M$, and
$\T^*(M(s))=(\T^* M)(s)$. The cotangent bundle $\T^*M(s)$ is a symplectic submanifold of
$\T^*M$ and the restriction $\omega\vert_{\mathbf{T}^*M(s)}$ is equal to the
Liouville one form $\omega_{s}$ on $\T^*M(s)$. The manifolds $\T^*M(s)$ are oriented
by their symplectic form $d\omega_{s}$.

For any $s\in K$, the tangent bundle $\T M$, when restricted to $M(s)$, decomposes as
$$
\T M\vert_{M(s)}=\T M(s)\oplus \Ncal.
$$
Let $\underline{s}$ be the linear action induced by $s$ on the bundle $\T M\vert_{M(s)}$: here
$\T M(s)$ is the kernel of $\underline{s}-{\rm Id}$, and the normal bundle $\Ncal$ is equal to the
image of $\underline{s}-{\rm Id}$.

Let $\nabla$ be a $K$-equivariant  connection on the the tangent bundle $\T M$. It induces
$K(s)$-equivariant connections : $\nabla^{0,s}$ on the bundle $\T M(s)$ and $\nabla^{1,s}$
on the bundle $\Ncal$. For $i=0,1$, we consider the equivariant curvature $R_i(Y), Y\in\kgot(s)$
of the connections $\nabla^{i,s}$.  We will use the following equivariant forms

\begin{defi}\label{def:JD-form}
We consider the following smooth closed $K(s)$-equivariant forms on $M(s)$:
\begin{eqnarray*}
\Ag(M(s))^2(Y)&:=& \det\left(\frac{R_0(Y)}{\e^{R_0(Y)/2}-\e^{-R_0(Y)/2}}\right), \\
\Df_s(\Ncal)(Y)&:=& \det\left(1-\underline{s}\e^{R_1(Y)}\right).
\end{eqnarray*}
The form  $\Ag(M(s))^2(Y)$  is defined for $Y$ in a (small)
neighborhood $\Ucal_s$ of $0\in\kgot(s)$.
\end{defi}

If the manifold $M(s)$ admits a $K(s)$-invariant orientation, one can define the square root
of $\Ag(M(s))^2$ : it is the equivariant $\Ag$-genus of the manifold $M(s)$.

The manifold $M(s)$ has several connected components $C_i$. We denote
by $\dim M(s)$ the locally constant function on $M(s)$ equal to $\dim C_i$ on $C_i$.
In the formulas of the cohomological index, we will use the following closed equivariant form on $M(s)$.

\begin{defi}\label{def:Lambda-s}
We consider the smooth closed equivariant form on $M(s)$
$$
\Lambda_s(Y):= (2i\pi)^{-\dim M(s)} \,  \frac{\Ag(M(s))^2(Y)}{\Df_s(\Ncal)(Y)}
$$
which is defined for $Y$ in a (small) neighborhood $\Ucal_s$ of $0\in\kgot(s)$.
\end{defi}

Here $\Ucal_s$ is a small $K(s)$-invariant neighborhood of $0$ in $\kgot(s)$.
It is chosen so that : $\ad(k)\Ucal_s=\Ucal_{ksk^{-1}}$,  and $M(s)\cap M(S)=M(s\e^S)$
for any $s\in K$ and any $S\in\Ucal_s$. For any $s\in K$ and any $S\in \Ucal_s$, let
$\Ncal_{(s,S)}$
be the normal bundle of $M(s\e^S)=M(s)\cap M(S)$ in $M(s)$. Let $R(Z)$ be the
$K(s\e^S)$-equivariant curvature of an invariant euclidean connection on $\Ncal_{(s,S)}$. Let
$$
\eul\left(\Ncal_{(s,S)}\right)(Z):=(-2\pi)^{-{\rm rank}\,\Ncal_{(s,S)}/2}\ \det{}^{1/2}_{o}( R(Z))
$$
be its $K(s\e^S)$-equivariant Euler form. Recall that $S$ induces a complex structure
on the bundle $\Ncal_{(s,S)}$. The square root $\det{}^{1/2}_{o}$ is computed using
the orientation $o$ defined by this complex structure.

Note that the diffeomorphism $k:\T^*M(s)\to\T^*M(ksk^{-1})$ induces a map
$k:\Acal^{\infty}_c(\Ucal_s, \T^*M(s))\to\Acal^{\infty}_c(\Ucal_{ksk^{-1}}, \T^*M(ksk^{-1}))$.
It is easy to check that the family $\Lambda_s\in\Acal^{\infty}(\Ucal_s, M(s))$ satisfies :
\begin{equation}\label{eq:forme-|ambda-1}
k\cdot\Lambda_s=\Lambda_{ksk^{-1}}\quad
\mathrm{in}\quad \Hcal^{\infty}(\Ucal_s, M(s)),
\end{equation}
\begin{equation}\label{eq:forme-lambda-2}
\Lambda_{s\e^S}(Z)=(-1)^r\frac{\Lambda_s\vert_{M(s\e^S)}}{\eul(\Ncal_{(s,S)})^2}(S+Z)\quad
\mathrm{in}\quad \Hcal^{\infty}(\Ucal', M(s\e^S)),
\end{equation}
where $\Ucal'\subset\kgot(s\e^S)$ is a small invariant neighborhood of $0$, and
$r=\frac{1}{2}{\rm rank}_\Rbb \,\Ncal_{(s,S)}$.

\bigskip

Let $\gamma_s\in \Acal^{\infty}_c(\Ucal(s),\T^* M(s))$ be a family of closed equivariant forms
with compact support. We look now at the family of {\em smooth} invariant functions
$$
\theta(\gamma)_s(Y)=\int_{\mathbf{T}^*M(s)}\Lambda_s(Y)\gamma_s(Y),\quad Y\in\Ucal_s.
$$
\begin{lem}\label{lem:Theta-gamma}
The family $\theta(\gamma)_s$ defines an invariant function $\Theta(\gamma)\in
C^{\infty}(K)^{K}$ if
\begin{equation*}\label{eq:bouquet-1}
k\cdot\gamma_s=\gamma_{ksk^{-1}}\quad
\mathrm{in}\quad \Hcal^{\infty}_c(\Ucal_s, \T^*M(s)),
\end{equation*}
and
\begin{equation*}\label{eq:bouquet-2}
\gamma_{s\e^S}(Z)=\gamma_s\vert_{\mathbf{T}^* M(s\e^S)}(S+Z)\quad
\mathrm{in}\quad \Hcal^{\infty}_c(\Ucal', \T^*M(s\e^S)),
\end{equation*}
where $\Ucal'\subset\kgot(s\e^S)$ is a small invariant neighborhood of $0$.
\end{lem}

\begin{proof}
The proof, that can be found in \cite{Duflo-Vergne-Asterisque} and
\cite{BV.inventiones.96.1,BV.inventiones.96.2},
follows directly from the localization formula in equivariant cohomology. Note that the square
$\eul(\Ncal_{(s,S)})^2$ is equal to the equivariant Euler form of the normal bundle of
$\T^*M(s\e^S)$ in $\T^*M(s)$.
\end{proof}

\medskip

In this article, the equivariant forms $\gamma_s$ that we use are the Chern forms
attached to a transversally elliptic symbol : since they have {\em generalized coefficients}, we will
need an extension of Lemma \ref{lem:Theta-gamma} in this case.

\subsection{The Chern character with support}\label{sec:chern-character-support}

Let $M$ be a $K$-manifold. Let $p:\T^*M\to M$ be the projection.

Let $\Ecal=\Ecal^+\oplus \Ecal^-$ be  a Hermitian $K$-equivariant super-vector bundle over $M$.
Let $\sigma:p^*\Ecal^+\to p^*\Ecal^-$ be a $K$-equivariant elliptic symbol: in this section
we do not impose any conditions of ellipticity on $\sigma$. Recall that $\supp(\sigma)\subset \T^*M$
is the set where $\sigma$ is not invertible.

Choose a $K$-invariant super-connection $\A$ on $p^*\Ecal$, without $0$ exterior degree
term. As in \cite{Quillen85,pep-vergneII}, we deform $\A$ with the help of
$\sigma$ : we consider the family of super-connections
$$
\A^{\sigma}(t)=\A+it\,v_\sigma,\,\, t\in \Rbb,
$$
on $\Ecal$ where
$v_\sigma=
\left(\begin{array}{cc}
0 & \sigma^*\\
\sigma & 0\\
\end{array}\right)$
is an odd endomorphism of $\Ecal$ defined with the help of the Hermitian structure. Let
$\F(\sigma,\A,t)(X), X\in \kgot$, be the equivariant curvature of
$\A^{\sigma}(t)$.

We denote by $\F(X),X\in \kgot$, the equivariant curvature
of $\A$: we have $\F(X)=\A^2+ \mu^\A(X)$ where $\mu^{\A}(X)\in \Acal(\T^*M,\End(p^*\Ecal))$ is
the moment of $\A$ \cite{BGV}. Then  $\F(\sigma,\A,t)(X)\in
\Acal(\T^*M,\End(p^*\Ecal))^+$ is given by:
\begin{equation}\label{F-sigma}
\F(\sigma,\omega,\A,t)(X)= -t^2v_\sigma +it[\A, v_\sigma]+\F(X).
\end{equation}

Let $t\in \Rbb$. Consider the $K$-equivariant forms on $\T^*M$:
\begin{eqnarray*}
\ch(\A)(X)&=&\str(\e^{\F(X)}),\\
\ch(\A,t)(X)&=& \str\left(\e^{\F(\sigma,\A,t)(X)}\right),\\
\eta(\sigma,\A,t)(X)&=&-i\str\left(v_\sigma \e^{\F(\sigma,\A,t)(X)}\right),\\
\beta(\sigma,\A,t)(X)&=&\int_{0}^{t}\eta(\sigma,\A,r)(X)dr.
\end{eqnarray*}
The forms $\ch(\A)$, $\ch(\A,t)$ and
$\beta(\sigma,\A,t)$ are equivariant forms on $\T^*M$ with
$C^{\infty}$-coefficients.  We have on $\T^* M$ the relation
$D(\beta(\sigma,\A,t))=\ch(\A)-\ch(\A,t)$.

We show in \cite{pep-vergneII} that the equivariant forms $\ch(\A,t)$ and $\eta(\sigma,\A,t)$ tends to zero
exponentially fast on the open subset $\T^*M\setminus\supp(\sigma)$, when $t$ goes to infinity.
Hence the integral
$$
\beta(\sigma,\A)(X)=\int_{0}^{\infty}\eta(\sigma,\A,t)(X)dt
$$
defines an equivariant form with $\f$-coefficients on $\T^*M\setminus\supp(\sigma)$, and we have
$D\beta(\sigma,\A)=\ch(\A)$ on $\T^*M\setminus \supp(\sigma)$.

We will now define the Chern character with support of $\sigma$.
For any invariant open neighborhood $U$ of $\supp(\sigma)$, we
consider the algebra $\Acal_U(\T^*M)$ of differential forms on
$\T^*M$ which are supported in $U$. Let
$\Acal_U^\infty(\kgot,\T^*M)$ be the vector space of equivariant
differential forms $\alpha:\kgot\to\Acal_U(\T^*M)$ which are
supported in $U$ :  $\Acal_U^\infty(\kgot,\T^*M)$ is a subspace of
$\Acal^\infty(\kgot,\T^*M)$ which is stable under the derivative
$D$. Let $\Hcal_U^\infty(\kgot,\T^*M)$ be the corresponding
cohomology space.

The following proposition follows
easily:

\begin{prop}\label{prop:ch-supp}
Let $U$ be a $K$-invariant open neighborhood of $\supp(\sigma)$.
Let $\chi\in\f(\T^* M)$ be a $K$-invariant function, with
support contained in $U$ and equal to $1$ in a neighborhood of
$\supp(\sigma)$. The equivariant differential form on $\T^*M$
\begin{eqnarray*}
  c(\sigma,\A,\chi)=\chi \ch(\A)+
   d\chi\,\beta(\sigma,\A)
\end{eqnarray*}
is equivariantly closed and supported in $U$. Its cohomology class $\ch_U(\sigma)$ in \break
$\Hcal^{\infty}_{U}(\kgot,\T^*M)$ does not depend on the choice of $(\A,\chi)$, nor on the Hermitian
structure on $\Ecal$.
\end{prop}

\begin{defi} We define the ``Chern character with support" $\chs(\sigma)$ as the collection
$(\ch_U(\sigma))_U$, where $U$ runs over $K$-invariant open neighborhood of $\supp(\sigma)$.
\end{defi}

In practice,  the Chern character with support $\chs(\sigma)$ will be identified with a class
$\ch_U(\sigma)\in \Hcal^{\infty}_{U}(\kgot,\T^*M)$, where $U$ is a ``sufficiently" small neighborhood
of $\supp(\sigma)$.

When $\sigma$ is elliptic,
we can choose $\chi\in\f(\T^*M)^K$ with compact support, and we denote
\begin{equation}\label{eq:ch-c-elliptic}
\ch_c(\sigma)\in \Hcal^{\infty}_c(\kgot,\T^*M)
\end{equation}
the class defined by the equivariant form with {\em compact support} $c(\sigma,\A,\chi)$.

\medskip

We introduce now the Bouquet of Chern characters with support.

Let $s\in K$. Then the action of $s$ on $\Ecal|_{M(s)}$ is given
by $s^{\Ecal}$, an even endomorphism of $\Ecal|_{M(s)}$. The
restriction of $\omega$ to $\T^*M(s)$ is the canonical $1$-form
$\omega_s$ of $\T^*M(s)$.

The super-connection $\A+itv_\sigma$ restricts to a
super-connection on $p^*\Ecal|\T^*M(s)$. Its curvature
$\F(\sigma,\A,t)$ restricted to $\T^*M(s)=N(s)$  gives an
element of \break $\Acal(N(s),\End(p^*\Ecal|_{N(s)}))$. To avoid further
notations, if $\chi$ is a function on $M$, we still denote by
$\chi$ its restriction to $M(s)$,  by $\sigma$ the restriction of
$\sigma$ to $\T^*M(s)$, by $\F(\sigma,\A,t)$ the
restriction of $\F(\sigma,\A,t)$ to $\T^*M(s)$, etc....

For $Y\in \kgot(s)$, we introduce the following $K(s)$-equivariant
forms on $\T^*M(s)$:

\begin{eqnarray*}
\ch_s(\A)(Y)&=&\str(s^{\Ecal}\e^{\F(Y)}),\\
\eta_s(\sigma,\A,t)(Y)&=&-i
   \str\left(v_\sigma\, s^{\Ecal}\,\e^{\F(\sigma,\A,t)(Y)}\right),\\
\beta_s(\sigma,\A)(Y)&=&\int_{0}^{\infty}\eta_s(\sigma,,\A,t)(Y)dt.
  \end{eqnarray*}
Then $\beta_s(\sigma,\omega,\A)$ is well defined $K(s)$-equivariant form with $\f$-coefficients
on  $\T^*M(s)\setminus\supp(\sigma)\cap \T^*M(s)$. We have similarly
$d\beta_s(\sigma,\A)=\ch_s(\A)$ outside $\supp(\sigma)\cap
\T^*M(s)$.

 The bouquet of  Chern characters $(\chs(\sigma,s))_{s\in K}$
can be constructed as follows.

\begin{prop}\label{prop:ch-supp-s}
Let $U$ be a $K(s)$-invariant open neighborhood of $\supp(\sigma)\cap \T^*M(s)$ in
$\T^*M(s)$. Let $\chi\in\f(\T^* M(s))$ be a $K(s)$-invariant function, with
support contained in $U$ and equal to $1$ in a neighborhood of
$\supp(\sigma)\cap \T^*M(s)$. The equivariant differential form on $\T^*M(s)$
\begin{eqnarray*}
  c_s(\sigma,\A,\chi)(Y)=\chi \ch_s(\A)(Y)+
   d\chi\,\beta_s(\sigma,\A)(Y),\quad Y\in\kgot(s),
\end{eqnarray*}
is equivariantly closed and supported in $U$. Its cohomology class $\ch_{U}(\sigma,s)$ in \break
$\Hcal^{\infty}_{U}(\kgot,\T^*M(s))$ does not depend on the choice of $(\A,\chi)$, nor on the Hermitian
structure on $\Ecal$.

We define the ``Chern character with support" $\chs(\sigma,s)$ as the collection
$(\ch_{U}(\sigma,s))_U$, where $U$ runs over the $K(s)$-invariant open neighborhood
of $\supp(\sigma)\cap \T^*M(s)$ in $\T^*M(s)$.
\end{prop}

\bigskip

\begin{lem}\label{resc}
Let $s\in K$ and $S\in K(s)$. Then for all $Y\in \kgot(s)\cap
\kgot(S)$, one has
$$
c_{s\e^S}(\sigma,\A,\chi)(Y)=c_{s}(\sigma,\A,\chi)(S+Y)|_{N(s)\cap N(S)}.
$$
\end{lem}

\begin{proof}
Let $N=\T^*M$. We have to compare the following forms on $N(s)\cap
N(S)$
\begin{eqnarray*}
\ch_{s\e^S}(\A)(Y)&=&\str(s^{\Ecal}\e^{S^{\Ecal}}\e^{\F(Y)}),\\
\ch_{s}(\A)(S+Y)&=&\str(s^{\Ecal}\e^{\F(S+Y)}),
\end{eqnarray*}
as well as the following forms
\begin{eqnarray*} \eta_{s\e^S}(\sigma,\omega,\A,t)(Y)&=& -i
   \str\left(v_\sigma \, s^{\Ecal}\e^{S^{\Ecal}}\e^{\F(\sigma,\A,t)(Y)}\right),\\
\eta_{s}(\sigma,\A,t)(S+Y)&=& -i
   \str\left(v_\sigma \, s^{\Ecal}\e^{\F(\sigma,\A,t)(S+Y)}\right).\\
\end{eqnarray*}

For $S\in \kgot$, the equivariant curvature
$\F(\sigma,\A,t)(S+Y)$ on $N(S)$ is equal to
$S^{\Ecal}+\F(\sigma,\A,t)(Y)$ as the vector field $VS$
vanishes on $N(S)$. Furthermore, above $N(s)\cap N(S)$, the
endomorphism  $\F(\sigma,\A,t)(Y)$ commutes with
$S^{\Ecal}$, for $Y\in \kgot(S)\cap \kgot(s)$. Thus the result
follows.
\end{proof}

\medskip

Then, for any open neighborhood $U$ of $\supp(\sigma)$, the family
$(\ch_{U}(\sigma,s))_{s\in K}$ forms a bouquet of cohomology
classes in the sense of \cite{Duflo-Vergne-Asterisque}.

\subsection{The Chern character of a transversally elliptic symbol}\label{sec:chern-character}

We keep the same notations than in the previous sections.
We denote by   $\omega$ the Liouville form  on $\T^* M$.
In local coordinates $(q,p)$ then $\omega=-\sum_{a}p_a dq_a$. The
two-form $\Omega=d\omega=\sum_a dq_a\wedge dp_a$ gives a
symplectic structure to $\T^*M$. The orientation of $\T^*M$ is the
orientation determined by the symplectic structure (our convention
for the canonical one-form $\omega$ differs from
\cite{BV.inventiones.96.1}, but the symplectic form $\Omega$ is
the same).

The moment map for the action of $K$ on $(\T^*M,\Omega)$ is the map
$f_\omega:\T^*M\to \kgot^*$ defined by
$\langle f_\omega(x,\xi), X\rangle=\langle \xi, V_xX\rangle$: we have
$D\omega(X)=\Omega + \langle f_\omega, X\rangle$.

Remark that  $\T^*_K M$ is the set of zeroes of
$f_\omega$. Recall how to associate to the $1$-form $\omega$ a
$K$-equivariant form $\Par(\omega)$ with generalized coefficients
supported near $\T_{K}^*M$.

On the complement of $\T_K^*M$, the $K$-equivariant form
\begin{equation}\label{eq:beta-omega}
\beta(\omega)=-i\omega\int_0^{\infty}\e^{it D\omega}dt
\end{equation} is well
defined as a $K$-equivariant form  with generalized coefficients,
and it is obvious to check that $D\beta(\omega)=1$ outside
$\T_{K}^*M$.

\begin{defi}\label{def:par-omega}
Let $U'$ be a $K$-invariant open neighborhood of $\T^*_KM$.
Let $\chi'\in\f(\T^* M)$ be a $K$-invariant function, with
support contained in $U'$ and equal to $1$ in a neighborhood of
$\T^*_KM$. The equivariant differential form on $\T^*M$
\begin{eqnarray*}
  \Par(\omega,\chi')=\chi'+
   d\chi'\,\beta(\omega)
\end{eqnarray*}
is closed, with generalized coefficients, and supported in $U'$. Its cohomology class
$\Par_{U'}(\omega)$ in $\Hcal^{-\infty}_{U'}(\kgot,\T^*M)$ does not depend on the choice of $\chi'$.

We will denote $\Par(\omega)$ the collection $(\Par_{U'}(\omega))_{U'}$.
\end{defi}

It is immediate to verify that
\begin{equation}\label{Parequalone}
\Par(\omega,\chi')=1+D\Big((\chi'-1)\beta(\omega)\Big).
\end{equation}

Thus, if we do not impose support conditions, the $K$-equivariant
form $\Par(\omega,\chi')$ represents $1$ in
$\Hcal^{-\infty}(\kgot,\T^*M)$.

\bigskip

We consider now a $K$-transversally elliptic symbol $\sigma$ on $M$.
We have the Chern character $\chs(\sigma)$ which is an equivariant form with
$\f$-coefficients which is supported near $\supp(\sigma)$, and the equivariant
form $\Par(\omega)$ with $\fgene$-coefficients which is supported near $\T^*_K M$.
Since $\supp(\sigma)\cap \T^*_K M$ is compact, the product
$$
\chs(\sigma)\wedge \Par(\omega)
$$
defines an equivariant form {\bf with compact support} with $\fgene$-coefficients.
We summarize the preceding discussion by the

\begin{theo}
Let $\sigma$ be a $K$-transversally elliptic symbol. Let $U,U'$ be respectively
$K$-invariant open neighborhoods of $\supp(\sigma)$ and $\T^*_K M$ such that
$\overline{U\cap U'}$ is compact. The product
$$
\ch_U(\sigma)\wedge\Par_{U'}(\omega)
$$
defines a compactly supported class in $\Hcal^{-\infty}_c(\kgot,\T^*M)$ which depends
uniquely of $[\sigma\vert_{\mathbf{T}^*_K M}]\in \KK_K(\T^*_K M)$ : this
equivariant class is denoted $\chc(\sigma,\omega)$.
\end{theo}

We will use the notation
$$
\ch_c(\sigma,\omega)=\chs(\sigma)\wedge\Par(\omega)
$$
which summarizes the fact that the class with compact support
$\chc(\sigma,\omega)$ is represented by the product
\begin{equation}\label{eq:rep-ch-c-transversal}
c(\sigma,\A,\chi)\wedge \Par(\omega,\chi')
\end{equation}
where $\chi,\chi'\in\f(\T^*M)^K$ are equal to $1$ respectively in a neighborhood of
$\supp(\sigma)$ and $\T^*_K M$, and furthermore the product $\chi\chi'$ is compactly supported.

\medskip

\begin{rem}
If $\sigma$ is elliptic, one can take $\chi$ with compact support, and $\chi'=1$ on $\T^*M$ in
Equation (\ref{eq:rep-ch-c-transversal}). We see then that
$$
\chc(\sigma,\omega)=\chc(\sigma)\quad \mathrm{in}\quad \Hcal^{-\infty}_c(\kgot,\T^*M).
$$
\end{rem}

\bigskip

Let $s\in K$. Similarly, we denote by $\Par(\omega_s,\chi')$ the
closed $K(s)$-equivariant form  on $\T^*M(s)$  associated to the
canonical $1$-form $\omega_s=\omega|\T^{*}M(s)$ and a function
$\chi'\in\f(\T^*M(s))^{K(s)}$ equal to $1$ in  a neighborhood of $\T^*_{K(s)}M(s)$.
For any $K(s)$-invariant neighborhood $U'\subset \T^*M(s)$ of $\T^*_{K(s)}M(s)$, we denote
$$
\Par_{U'}(\omega_s)\in\Hcal^{-\infty}_{U'}(\kgot(s),\T^*M(s))
$$
the class defined by $\Par(\omega_s,\chi')$ when $\chi'$ is supported in $U'$.
We denote $\Par(\omega_s)$ the collection $(\Par_{U'}(\omega,s))_{U'}$.

We defined, in Section \ref{sec:chern-character-support}, the
family of Chern classes $(\chs(\sigma,s))_{s\in K}$ for any
$K$-invariant symbol. We now define a family
$(\chc(\sigma,\omega,s))_{s\in K}$ with compact support and
$\fgene$-coefficients when $\sigma$ is a $K$-transversally
elliptic symbol on $M$. Note that the restriction
$\sigma|\T^{*}M(s)$ of a $K$-transversally elliptic symbol on $M$
is a $K(s)$-transversally elliptic symbol on $M(s)$.

\begin{prop}\label{prop:ch-c-omega-s}
Let $\sigma$ be a $K$-transversally elliptic symbol. Let $s\in K$.

Let $U,U'\subset \T^*M(s)$ be respectively
$K(s)$-invariant open neighborhoods of $\supp(\sigma|\T^{*}M(s))$ and $\T^*_{K(s)} M$ such that
$\overline{U\cap U'}$ is compact. The product
$$
\ch_U(\sigma,s)\wedge\Par_{U'}(\omega_s)
$$
defines a compactly supported class in $\Hcal^{-\infty}_c(\kgot(s),\T^*M(s))$ which depends
uniquely of $[\sigma\vert_{\mathbf{T}^*_{K(s)} M(s)}]$. The notation
$$
\ch_c(\sigma,\omega,s)=\chs(\sigma,s)\wedge\Par(\omega_s)
$$
summarizes the fact that the class with compact support
$\chc(\sigma,\omega,s)$ is represented by $c_s(\sigma,\A,\chi)\wedge \Par(\omega_s,\chi')$
where $\chi,\chi'\in\f(\T^*M(s))^{K(s)}$ are chosen so that $\chi\chi'$ is compactly supported.
\end{prop}

\subsection{Definition of the cohomological index}

Let $K$ be a compact Lie group and let $M$ be a compact $K$-manifold. The
aim of this section is to define the cohomological index
$$
\index^{K,M}_c: \KK_K(\T^*_K M)\to\fgene(K)^{K}.
$$

For any $[\sigma]\in\KK_K(\T^*_K M)$, the generalized function $\index^{K,M}_c([\sigma])$ will
be described through their restrictions $\index^{K,M}_c([\sigma])\vert_s, s\in K$ (see Section
\ref{sec:restriction}).

In Subsection \ref{subsec:integration-bouquet}, we have introduce for any $s\in K$,
the closed equivariant form on $M(s)$
$$
\Lambda_s(Y):= (2i\pi)^{-\dim M(s)} \,  \frac{\Ag(M(s))^2(Y)}{\Df_s(\Ncal)(Y)}.
$$

We wish to prove first the following theorem.

\begin{theo}\label{likeAS}
Let $\sigma$ be a $K$-transversally elliptic symbol. There exists a unique
invariant generalized function $\index_c^{K,M}(\sigma)$
on $K$ satisfying the following equations.
 Let $s\in K$. For every  $Y\in \kgot(s)$ sufficiently small,
\begin{equation}\label{spart}
\index_c^{K,M}(\sigma)\vert_s(Y)=\int_{\mathbf{T}^*M(s)}\Lambda_s(Y)\,\chs(\sigma,s)(Y)\,
\Par(\omega_s)(Y).
\end{equation}
\end{theo}

As $\chc(\sigma,\omega,s)=\chs(\sigma,s)(Y)\, \Par(\omega_s)(Y)$
is {\bf compactly supported}, the integral (\ref{spart}) of
equivariant differential forms with generalized coefficients
defines a generalized function on a neighborhood of zero in
$\kgot(s)$. However, we need to prove that the different local
formulae match together. The proof of this theorem occupies the
rest of this subsection. Once this theorem is proved, we can make
the following definition.

\begin{defi}\label{likeASdefi}
Let $\sigma$ be a $K$-transversally elliptic symbol. The
cohomological index of $\sigma$ is the  invariant generalized
function $\index_c^{K,M}(\sigma)$ on $K$ satisfying Equation
(\ref{spart}). We also rewrite the formula for the cohomological
index as
\begin{equation}\label{spart-bis}
\index_c^{K,M}(\sigma)\vert_s(Y)=
\int_{\mathbf{T}^*M(s)}\Lambda_s(Y)\,\chc(\sigma,\omega,s)(Y).
\end{equation}

In particular, when $s=e$, Equation (\ref{spart}) becomes
\begin{equation}\label{spart-e}
\index_c^{K,M}(\sigma)(e^X)=(2i\pi)^{-\dim M}\int_{\mathbf{T}^*M}
\Ag(M)^2(X)\chs(\sigma)(X)\,\Par(\omega)(X).
\end{equation}
\end{defi}
\begin{rem} In (\ref{spart}), (\ref{spart-e}) and (\ref{spart-bis}) we take
for the integration the symplectic orientation on the
cotangent bundles.
\end{rem}

\medskip

Let us now prove Theorem \ref{likeAS}.

\begin{proof}  The right hand side of (\ref{spart}) defines a
$K(s)$-invariant generalized function $\theta_s(Y)$ on a neighborhood
$\Ucal_s$ of $0$ in $\kgot(s)$. Following Theorem
\ref{theo:germe-function}, the family $(\theta_s)_{s\in K}$
defines an invariant generalized function on $K$, if the {\em invariance condition} and the {\em
compatibility condition} are satisfied. The invariance condition
is easy to check. We will now prove the compatibility condition.

Let $s\in K$, and $S\in\Ucal_s$. We have to check that the
restriction $\theta_s\vert_S$ coincides with $\theta_{s\e^S}$ in a
neighborhood of $0$ in $\kgot(s)\cap\kgot(S)=\kgot(s\e^S)$. We
conduct the proof only for $s$ equal to the identity $e$, as the
proof for $s$ general is entirely similar.

We have to compute the restriction at $\theta_e\vert_S$ of the
generalized invariant function
$$
\theta_e(X):=\int_{\mathbf{T}^*M}\Lambda_e(X)\,\chs(\sigma)(X)\,
\Par(\omega)(X), \quad X\in\Ucal_e.
$$
For this purpose, we choose a particular representant of the class
$\chs(\sigma)\Par(\omega)$ in $\Hcal^{-\infty}_c(\kgot,\T^*M)$.
Since this class only depends of
$[\sigma\vert_{\mathbf{T}^*_K M}]\in\KK_K(\T^*_K M)$, we choose a
transversally elliptic symbol $\sigma_h$ which is almost
homogeneous of degree $0$ and such that
$[\sigma_h\vert_{\mathbf{T}^*_K M}]= [\sigma\vert_{\mathbf{T}^*_K
M}]$.

One can show \cite{pep-vergneIII} that the moment map
$f_\omega:\T^*M\to\kgot^*$ is {\em proper} when restricted to the
support $\supp(\sigma_h)$.  We represent $\chs(\sigma_h)$ by the form
$c(\sigma_h,\A,\chi_\sigma)$ where  $\chi_\sigma$ is a function on $\T^*M$
such that ${\rm support}(\chi_\sigma)\cap \{\|f_\omega\|^2\leq  1\}$ is compact.
For this choice of $\chi_\sigma$,
the equivariant form $\alpha(X)=\Lambda_e(X)\,c(\sigma_h,\A,\chi_\sigma)(X)$ is
thus such that ${\rm support}(\alpha)\cap \{\|f_\omega\|^2\leq 1\}$ is
compact. It is defined for $X$ small enough. Multiplying by a
smooth invariant function of $X$ with small compact support and
equal to $1$ in a neighborhood of $0$, we may find  $\alpha(X)$
defined for all $X\in \kgot$ and which coincide with
$\Lambda_e(X)\,c(\sigma_h,\A,\chi_\sigma)(X)$ for $X$ small enough.

We choose $\chi$ supported in $\{\|f_\omega\|^2< 1\}$ and  equal to $1$
on $\{\|f_\omega\|^2\leq \epsilon\}$, and define $\Par(\omega)$ with this choice of
$\chi$. Then $\alpha(X) \Par(\omega)(X)$ is compactly supported.

We will now prove the following  result:

\begin{prop}\label{prop-loc-omega}
Let $\alpha(X)$ be a closed equivariant form with
$C^{\infty}$-coefficients on $N=\T^*M$ such that $\{\|f_\omega\|^2\leq 1\}\cap{\rm support}(\alpha)$
is compact. Define the generalized
function $\theta\in C^{-\infty}(\kgot)^K$ by
\begin{equation}\label{eq:prop-loc-omega-1}
\theta(X):=\int_{N}\alpha(X) \Par(\omega)(X).
\end{equation}

 Then, the restriction $\theta|_S$ is given, for $Y=Z-S$ sufficiently close to $0$ by
\begin{equation}\label{eq:prop-loc-omega-2}
\theta|_S(Y)=(-1)^{r}\int_{N(S)}\frac{\alpha\vert_{N(S)}(Z)}{\eul(\Ncal_{S})^2(Z)}
\Par(\omega_S)(Z).
\end{equation}

Here $\Ncal_S$ denotes the normal bundle of $M(S)$ in $M$, and $r=\frac{1}{2}(\dim M-\dim M(S))$.

\end{prop}

\begin{rem} The integral (\ref{eq:prop-loc-omega-2}) is defined using
the symplectic orientation $o(\omega_S)$ on
$N(S)=\T^*M(S)$. The linear action of $S$ on the normal bundle $\Ncal'_S$ of $N(S)$ in $N$ induces a
complex structure $J_S$ :  let $o(J_S)$ be corresponding orientation of the fibers of $\Ncal'_S$.
We have then on $N(S)$ the orientation $o(S)$ such that $o(\omega)=o(S)o(J_S)$. One can check
that $(-1)^r$ is the quotient between $o(S)$ and $o(\omega_S)$.
\end{rem}

Let us apply the last proposition to the form $\alpha(X)=\Lambda_e(X)\chs(\sigma_h)(X)$. If we use
(\ref{eq:forme-lambda-2}) and Lemma \ref{resc}, we see that
$(-1)^{r}\frac{\alpha\vert_{N(S)}}{\eul(\Ncal_{S})^2}(S+Y)$
is equal to \break $\Lambda_{\e^S}(Y)\chs(\sigma_h,\e^S)(Y)$. Hence Proposition \ref{prop-loc-omega}
tells us that the the restriction of $\theta_e\vert_S$ is equal to $\theta_{\e^S}$ : Theorem \ref{likeAS}
is proved.

\medskip

\begin{proof} We now concentrate on the proof of Proposition \ref{prop-loc-omega}.

Remark that if $\alpha$ is compactly supported, we can get rid of the forms
$\Par(\omega)$ and $\Par(\omega_S)$ in the integrals (\ref{eq:prop-loc-omega-1}) and
(\ref{eq:prop-loc-omega-2}), since they are equal to $1$ in cohomology.  In this case, the proposition is just
the localization formula, as $\eul(\Ncal_{S})^2$ is the Euler
class of the normal bundle of $ \T^*M(S)$ in $\T^*M$.

The proof will follow the same scheme as the
usual localization formula (see \cite{BGV}) and will use  the fact that
$\alpha|_{\kgot(S)}$ is exact outside the set of zeroes of $S$. To
extend the proof  of the localization formula in our setting, we have to bypass the fact that
the restriction of $\Par(\omega)$ to $\kgot(S)$ has no meaning, since $\Par(\omega)$ is an equivariant form
with generalized coefficients. However, we will use in a crucial
way the fact that the closed equivariant form $\Par(\omega)$ is the limit of smooth
equivariant forms
$$
\Par^T(\omega)(X)=\chi+d\chi\int_0^T(-i\omega \e^{itD\omega(X)})dt.
$$
Here $D\left(\Par^T(\omega)\right)=d\chi\e^{iTD\omega}$ tends to zero
as $T$ goes to infinity.

Let $\Par^T(\omega)(Z)$ be the restriction of $\Par^T(\omega)(X)$ to $\kgot(S)$.  We write
$f_\omega=f^{S}_\omega+ f^\qgot_\omega$ relative to the $K(S)$-invariant decomposition
$\kgot^*=\kgot(S)^*\oplus\qgot$. Then the family of $K(S)$-equivariant forms
$$
\e^{itD\omega(Z)}=\e^{itd\omega}\e^{it\langle f_\omega^S, Z\rangle}
$$
tends to $0$ outside $\{f_\omega^S=0\}$, as $t$ goes to $\infty$.
Since $d\chi$ can be non-zero on the subset $\{f_\omega^S=0\}$, the family of $\kgot(S)$-equivariant
forms $\Par^T(\omega)(Z)$ does not have a limit when $T\to\infty$ in general.

Consider the sub-manifold $N(S):=\T^*M(S)$ of $N:=\T^*M$. Note that $f^\qgot_\omega$ vanishes on
$N(S)$. Let $\Vcal$ be an invariant tubular neighborhood of $N(S)$ which is contained in
$\{\|f^\qgot_\omega\|^2\leq \frac{\epsilon}{2}\}$. We are interested in the restriction
$\Par^T(\omega)\vert_\Vcal(Z)$ to $\Vcal$. Since the function $\chi$ is equal to $1$ on
$\{\|f_\omega\|^2\leq\epsilon \}$, we see that $d\chi\vert_\Vcal$ is equal to
zero in the neighborhood $\Vcal\cap\{\|f^S_\omega\|^2\leq \frac{\epsilon}{2}\}$ of $\Vcal\cap\{f^S_\omega=0\}$.
Hence the limit
\begin{equation}\label{eq:limit-par-Vcal}
\Par(\omega)\vert_\Vcal(Z)=\lim_{T\to\infty}\Par^T(\omega)\vert_\Vcal(Z),\quad Z\in\kgot(S),
\end{equation}
defines a $K(S)$-equivariant form with generalized coefficients on $\Vcal$. Note that the restriction
of $\Par(\omega)\vert_\Vcal$ to $N(S)\subset \Vcal$ is the $K(S)$-equivariant form
$\Par(\omega_S)$ associated to the Liouville $1$-form $\omega_S$ on $\T^*M(S)$.

\medskip

The generalized function $\theta\in\fgene(\kgot)^K$ is the limit,
as $T$ goes to infinity, of the family of smooth functions
$$
\theta^T(X):= \int_{N} \alpha(X)\Par^T(\omega)(X).
$$
Here the  equivariant forms $\alpha^T=\alpha\,\Par^T(\omega)$ stay supported  in the fixed compact
set $\Kcal:=\{\|f_\omega\|^2\leq 1\}\cap {\rm support}(\alpha)$.

The proof will be completed if we show that the family of smooth
functions $\theta^T(Z),Z\in \kgot(S)$,  converge to the
generalized function
$$
\theta'(Z):=(-1)^r\int_{N(S)}\frac{\alpha(Z)}{\eul(\Ncal_{S})^2(Z)}
  \Par(\omega_S)(Z),
$$
as $T$ goes to infinity, and when $Z$ varies in a small neighborhood of $S$ in
$\kgot(S)$.

Let $U$ be a relatively compact invariant neighborhood of $\Kcal$
in $N$. Let $\chi'\in\f(U)^{K(S)}$ be such that $\chi'$ is
supported in $\Vcal\cap U$, and $\chi'=1$ in a neighborhood of
$U(S)=N(S)\cap U$.  Here $\Vcal$ is a tubular neighborhood of
$N(S)$ satisfying the conditions for the existence of the limit
(\ref{eq:limit-par-Vcal}).

Choose a $K$-invariant metric $\langle-,-\rangle$ on $\T N$. Let
$\lambda$ be the $K(S)$-invariant
one form on $N$ defined by $\lambda=\langle VS,-\rangle$. Note that $D(\lambda)(S)=d\lambda- \|VS\|^2$
is invertible outside $N(S)$. One sees that
$$
\mathrm{P}_{\chi'}(Z)=\chi'+d\chi'\frac{\lambda}{D\lambda(Z)}
$$
is a $K(S)$-equivariant form on $U$ for $Z$ in a small neighborhood of $S$. The following equation
of $K(S)$-equivariant forms on $U$ is immediate to verify:
\begin{equation}\label{eq:Par=1}
1= \mathrm{P}_{\chi'} +D\left((1-\chi')\frac{\lambda}{D\lambda}\right).
\end{equation}
Since the $K(S)$-equivariant forms
$$\alpha^T(Z):=\alpha(Z)\Par^T(\omega)(Z)$$ are supported in $U$,
one can multiply (\ref{eq:Par=1}) by $\alpha^T$. We have then the
following relations between {\em compactly supported}
$K(S)$-equivariant forms on $N$:
\begin{eqnarray*}
\alpha^T&=& \mathrm{P}_{\chi'}\alpha^T +D\left((1-\chi')\frac{\lambda}{D\lambda}\right)\alpha^T\\
&=&\mathrm{P}_{\chi'}\alpha^T +
D\left((1-\chi')\frac{\lambda}{D\lambda}\alpha^T\right) +
(1-\chi')\frac{\lambda}{D\lambda}D(\alpha^T).
\end{eqnarray*}
According to this equation, we divide the function $\theta^T(Z)$ in two parts
$$
\theta^T(Z)= A^T(Z)+ B^T(Z),\quad \mathrm{for}\quad Z-S \quad \mathrm{small}
$$
with
$$
A^T(Z)=\int_N \mathrm{P}_{\chi'}(Z)\alpha^T(Z)
$$
and
$$
B^T(Z)=\int_N (1-\chi')\frac{\lambda}{D\lambda(Z)}D\alpha^T(Z)=
\int_N (1-\chi')\frac{\lambda}{D\lambda(Z)}\alpha(Z) d\chi\e^{iTD\omega(Z)}.
$$

\medskip

Let $p: \Vcal\to N(S)$ be the projection, and let $i:N(S)\to
\Vcal$ be the inclusion. Since the form $\mathrm{P}_{\chi'}(Z)$ is
supported in $\Vcal$, the family of smooth equivariant forms
$\mathrm{P}_{\chi'}(Z)\alpha(Z)\Par^T(\omega)(Z)$ converges to
$$
\mathrm{P}_{\chi'}(Z)\alpha(Z)\Par(\omega)\vert_\Vcal(Z)
$$
as $T$ goes to $\infty$, by our previous computation of the limit (\ref{eq:limit-par-Vcal}).
Hence the functions $A^T(Z)$ converge to
\begin{eqnarray*}
\int_\Vcal \mathrm{P}_{\chi'}(Z)\alpha(Z)\Par(\omega)\vert_\Vcal(Z)
&=&\int_\Vcal \mathrm{P}_{\chi'}(Z)\ p^*\!\!\circ \! i^*\left(\alpha\,\Par(\omega)\vert_\Vcal\right)(Z) \quad [1]\\
&=&\int_{N(S),o(S)} p_*(\mathrm{P}_{\chi'})(Z)\ \alpha\vert_{N(S)}(Z)\,\Par(\omega_S)(Z) \quad [2]\\
&=&\int_{N(S),o(S)}\frac{\alpha\vert_{N(S)}(Z)}{\eul(\Ncal_S)^2(Z)} \Par(\omega_S)(Z) \quad [3]\\
&=&(-1)^r \int_{N(S)}\frac{\alpha\vert_{N(S)}(Z)}{\eul(\Ncal_S)^2(Z)} \Par(\omega_S)(Z). \quad [4]
\end{eqnarray*}
Points $[1]$  and $[2]$ are due to the fact that $\alpha(Z)\Par(\omega)\vert_\Vcal(Z)$ is equal to \break
$p^*\!\!\circ \! i^*\left(\alpha\,\Par(\omega)\vert_\Vcal\right)(Z)$ in $\Hcal^{-\infty}(\kgot(S),\Vcal)$ and that
$\mathrm{P}_{\chi'}$ has a compact support relatively to the
fibers of $p$ (here $p_*$ denotes the integration along the fibers). For point $[3]$, we use then that
$p_*(\mathrm{P}_{\chi'})$ multiplied
by the Euler class\footnote{The Euler form of the vector bundle $\Vcal\to N(S)$ is
equal to the square of the Euler form of the normal bundle $\Ncal_S$ of  $M(S)$ in  $M$.} of $\Vcal$
is equal to the restriction of $\mathrm{P}_{\chi'}$ to $N(S)$, which is identically equal to $1$.  In $[4]$, we use the
symplectic orientation for the integration.

\medskip

Let us show that the integral $\int_{\kgot(S)}B^T(Z)\varphi(Z)dZ$ tends to $0$, as $T$ goes to infinity,
for any  $\varphi\in\f(\kgot(S))^{K(S)}$ supported in a small neighborhood of $S$.
As $\det_{\kgot/\kgot(S)}(Z)$ does not vanish when $Z-S$ remains small enough, it is enough to show that
$$
I(T):=\int_{N\times \kgot(S)} (1-\chi')\frac{\lambda}{D\lambda(Z)}D\alpha^T(Z)
\varphi(Z)\det{}_{\kgot/\kgot(S)}(Z)dZ
$$
tends to $0$, as $T$ goes to infinity. We have
$$
I(T):=\int_{N\times \kgot(S)} \e^{i T D\omega(Z)} \eta(Z)\det{}_{\kgot/\kgot(S)}(Z)dZ
$$
where $\eta(Z)= (\chi'-1)\frac{\lambda}{D\lambda(Z)}\alpha(Z)d\chi\varphi(Z)$
is a compactly supported $K(S)$-equivariant form on $N$ with $\f$-coefficients, which
is defined for all $Z\in \kgot(S)$. Furthermore we have $\eta(Z)=0$ for $Z$ outside a small neighborhood of
$S$ and
$$
\mathrm{support}(\eta)\cap \big\{ f_{\omega}=0\big\}=\emptyset.
$$

There exists a $K$-equivariant form $\Gamma:\kgot\to \Acal(N)$
such that $\Gamma(Z)=\eta(Z)$ for any $Z-S$ small in $\kgot(S)$. Indeed we define
$\Gamma(X)=k\cdot \eta(Z)$ for any choice of $k,Z$ such that
$k\cdot Z=X$.  Here $X$ varies in  a (small) neighborhood of $K\cdot
S$. As $\eta(Z)$ is zero when $Z$ is not near $S$,  the map $X\mapsto\Gamma(X)$ is
supported on a compact neighborhood of $K\cdot S$ in $\kgot$.  We see also that
\begin{equation}\label{eq:support-gamma}
\mathrm{support}(\Gamma)\cap \big\{ f_{\omega}=0\big\}=\emptyset.
\end{equation}

Condition (\ref{eq:support-gamma}) implies that the integral
$J(T):=\int_{\kgot\times N}\e^{i T D\omega(X)}\Gamma(X) dX$ goes to
$0$, as $T$ goes to infinity. But $I(T)=J(T)$. Indeed, write $X=k\cdot Z$ and apply Weyl
integration formula. We obtain
\begin{eqnarray*}
J(T)&=&\int_{\kgot(S)}\Big(\int_{K\times N}\e^{i T D\omega(k\cdot Z)}
\Gamma(k\cdot Z) dk\Big)\det{}_{\kgot/\kgot(S)}(Z) dZ \\
&=&\int_{\kgot(S)}\int_{K\times N}k\cdot\Big (\e^{i T D\omega(Z)}\eta(Z)\Big)dk\
\det{}_{\kgot/\kgot(S)}(Z) dZ.
\end{eqnarray*}

Integration on the $K$-manifold $N$ is invariant by
diffeomorphisms, thus
$$J(T)=\int_{\kgot(S)}\int_{N}\e^{i T D\omega(Z)}\eta(Z)
\det{}_{\kgot/\kgot(S)}(Z) dZ =I(T).
$$

We have shown that the family of smooth function $B^T(Z)$ goes to $0$, as $T$ goes to infinity.
The proof of Proposition \ref{prop-loc-omega} is then completed.
\end{proof}

\end{proof}

\bigskip

\bigskip

Let $H$ be a compact Lie group acting on $M$ and commuting with
the action of $K$. Then the space $\T^*_KM$ is provided with an
action of $K\times H$.

\begin{lem}\label{lem:index-smooth-H}
If $[\sigma]\in \KK_{K\times H}(\T^*_KM)$, then the cohomological index \break
$\index^{K,H,M}_c(\sigma)\in C^{-\infty}(K\times H)^{K\times H}$ is smooth relatively to $H$.
\end{lem}

\begin{proof}
We have to prove that for any $s=(s_1,s_2)\in K\times H$, the generalized function
$$
\index^{K,H,M}_c(\sigma)\vert_s(Y_1,Y_2)
$$
which is defined for $(Y_1,Y_2)$ in a neighborhood of $0$ in
$\kgot(s_1)\times\hgot(s_2)$, is smooth relatively to the
parameter $Y_2\in\hgot(s_2)$. We check it for $s=e$.

We have
\begin{equation}\label{eq:indice-K-H}
\index^{K,H,M}_c(\sigma)\vert_e(X,Y)=\int_{\mathbf{T}^*M}\Lambda_e(X,Y)\chs(\sigma)(X,Y)
\Par(\omega)(X,Y)
\end{equation}
for $(X,Y)\in\kgot\times\hgot$ in a neighborhood of $0$.
The equivariant class with compact support $\chs(\sigma)\Par(\omega)$ is represented
by the product $c(\sigma,\A,\chi)\Par(\omega,\chi')$  where $(\chi,\chi')$ is chosen so that
$\chi=1$ in a neighborhood of $\supp(\sigma)$, $\chi'=1$  in a neighborhood of $\T^*_{K\times H}M$, and
$\chi\chi'$ is compactly supported.

Since $\sigma$ is $K$-transversally elliptic, the set $\supp(\sigma)\cap\T^*_{K}M$ is compact.
Hence we can choose $(\chi,\chi')$ so that $\chi'=1$  in a neighborhood of $\T^*_{K}M$ and
$\chi\chi'$ is compactly supported.  It easy to check that the equivariant form
$\Par(\omega,\chi')(X,Y)$ is then smooth relatively to the parameter $Y\in \hgot$.
This show that the right hand side
of (\ref{eq:indice-K-H}) is smooth  relatively to the parameter $Y\in \hgot$.

\end{proof}

 \begin{rem}\label{rem:chc-smooth-H}
 We will denote $\chcf(\sigma,\omega)(X,Y)$ the $K\times H$-equivariant form defined by the product
 $c(\sigma,\A,\chi)\Par(\omega,\chi')$  where $(\chi,\chi')$ is chosen so that
$\chi=1$ in a neighborhood of $\supp(\sigma)$, $\chi'=1$  in a neighborhood of $\T^*_KM$, and
$\chi\chi'$ is compactly supported. The equivariant form $\chcf(\sigma,\omega)(X,Y)$ is compactly supported and is
smooth relatively to $Y\in\hgot$.
\end{rem}

\section{The cohomological index coincides with the analytic one}

In this section, we now prove that the cohomological index is
equal to the analytical index. The main difficulty in the proof of
this result in Berline-Vergne \cite{BV.inventiones.96.1,BV.inventiones.96.2}
was to prove that their formulae were defining generalized functions which, moreover, were
compatible with each other. The heart of this new proof is
the fact the Chern character with compact support is multiplicative. Thus we rely
heavily here on the results of \cite{pep-vergneIII}, so that the proof is now easy.

\begin{theo}
The analytic index of a transversally elliptic operator $P$ on a $K$-manifold $M$ is equal to
$\index^{K,M}_c([\sigma_P])$.
\end{theo}

To prove that the cohomological index is equal to the analytic
index, following Atiyah-Singer algorithm, we need only to verify that the
cohomological index satisfies  the properties that we listed of
the analytic index:

\begin{itemize}

\item Invariance by diffeomorphism : $\mathbf{Diff}$,

\item Functorial with respect to subgroups : $\mathbf{Morph}$,

\item Excision property,

\item Free action properties,

\item Multiplicative properties,

\item Normalization conditions $\mathbf{[N1]},\mathbf{[N2]}$ and $\mathbf{[N3]}$.

\end{itemize}

The invariance by diffeomorphism, the functoriality with respect to subgroups and the excision property are
obviously satisfied by  $\index^{K,M}_c$.

\subsection{Free action}

We now prove that the cohomological index satisfies the free action property.
We consider the setting of Subsection \ref{subsec:a-index-free}. The action of
$K$ on the bundle $\T^*_K P$ is free and the quotient $\T^*_K P/K$ admit a canonical
identification with $\T^* M$. Then we still denote by
$$
q:\T^*_K P\to \T^* M
$$
the quotient map by $K$: it is a $G$-equivariant map such that $q^{-1}(\T^*_G M)=\T^*_{K\times G} P$.

We choose a $G$-invariant connection $\theta$ for the principal fibration $q: P\to M$ of group $K$.
With the help of this connection, we have a direct sum decomposition
$$
\T^* P= \T^*_K P\oplus P\times\kgot^*.
$$
Let $\pi_1:\T^* P\to \T^*_K P$ and $\pi_2:\T^* P\to P\times\kgot^*$ be the projections on each factors.
Let
$$
Q:\T^*P\to \T^*M
$$
be the map $q\circ \pi_1$.

Let $\sigma$ be a $G$-transversally elliptic morphism on $\T^*M$. Its pull-back
$Q^*\sigma$ is then a $K\times G$-transversally elliptic morphism on $\T^*P$ : we have
$\supp(Q^*\sigma)=Q^{-1}(\supp(\sigma))$ and then $\supp(Q^*\sigma)\cap\T^*_{K\times G} P
=q^{-1}(\supp(\sigma)\cap\T^*_{G} M)$ is compact.

\begin{theo}
Let $P\to M$ be a principal fibration with a free right action of $K$, provided with a left action of $G$.
Consider a class $[\sigma]\in \KK_G(\T^*_G M)$ and its pull-back by $Q$ :
$[Q^*\sigma]\in \KK_{K\times G}(\T^*_{K\times G}P)$. Then we  have the equality of generalized functions :
for $(k,g)\in K\times G$
$$
\index^{K\times G,P}_c([Q^*\sigma])(k,g)=\sum_{\tau\in \hat{K}} \tr(k,\tau)
\index^{G,M}_c([\sigma_{\tau^*}])(g).
$$
\end{theo}

The rest of this section is devoted to the proof. We have to check that for any $(s,s')\in K\times G$
we have the following equality of generalized functions defined in a neighborhood
of $\kgot(s)\times\ggot(s')$ :
\begin{equation}\label{eq:preuve-s-s'}
\index^{K\times G,P}_c([Q^*\sigma])\vert_{(s,s')}(X,Y)=\sum_{\tau\in \hat{K}} \tr(s\e^X,\tau)
\index^{G,M}_c([\sigma_{\tau^*}])\vert_{s'}(Y).
\end{equation}

We conduct the proof of (\ref{eq:preuve-s-s'}) only for $(s,s')=(1,1)$. This proof can be adapted to the general case
by using the same arguments as Berline-Vergne \cite{BV.inventiones.96.2}.

\bigskip

First, we analyze the left hand side of (\ref{eq:preuve-s-s'}) at
$(s,s')=(1,1)$.

We consider the
$K\times G$-invariant one form $\nu=\langle \xi,\theta\rangle$ on $P\times \kgot^*$ :
here $\theta\in\Acal^1(P)\otimes\kgot$
is our connection form, and $\xi$ is the variable in $\kgot^*$. We have
\begin{equation}\label{eq:forme-nu}
D\nu(X,Y)=d\nu +\langle \xi,\mu(Y)- X\rangle,\quad X\in\kgot,\quad Y\in\ggot.
\end{equation}
where $\mu(Y)=-\theta(V Y)\in C^{\infty}(P)\otimes \kgot$.

We associate to $\nu$ the $K\times G$-equivariant form with {\em generalized coefficients}
$\beta(-\nu)(X,Y)=i\nu\int_0^{\infty}\e^{-it D\nu(X,Y)}dt$,  $(X,Y)\in\kgot\times\ggot$,
 which is defined on the open subset $P\times \kgot^*\setminus\{0\}$. One checks that
$\beta(-\nu)(X,Y)$ is smooth relatively to the variable $Y\in\ggot$. Let $\chi_{\kgot^*}\in\f(\kgot^*)^K$
be a function with compact support and equal to $1$ near $0$. Then
\begin{equation}\label{eq:par-nu}
\Par(-\nu)(X,Y):=\chi_{\kgot^*}+d\chi_{\kgot^*}\beta(-\nu)(X,Y)
\end{equation}
is a closed equivariant form on $P\times\kgot^*$, with compact support, and which is smooth relatively to the
variable $Y\in\ggot$.

Let $\sigma$ be a $G$-transversally elliptic morphism on $\T^*M$.
Its pull-back $Q^*\sigma$ is then a $K\times G$-transversally
elliptic morphism on $\T^*P$. Let $\omega_P$ and $\omega_M$ be the
Liouville $1$-forms on $\T^*P$ and $\T^*M$ respectively. We have
defined the equivariant Chern class with compact support
$\ch_c(\sigma,\omega_M)\in \Hcal^{-\infty}_c(\ggot,\T^*M)$ and
$\ch_c(Q^*\sigma,\omega_P)\in\Hcal^{-\infty}_c(\kgot\times\ggot,\T^*P)$.

\begin{prop}\label{prop:ch-c-P-M}
We have the following equality
$$
\ch_c(Q^*\sigma,\omega_P)(X,Y)=
Q^*\Big(\ch_c(\sigma,\omega_M)\Big)(Y)\wedge
\pi_2^*\Big(\Par(-\nu)\Big)(X,Y)
$$
in $\Hcal^{-\infty}_c(\kgot\times\ggot,\T^*P)$. Note that the product on the right hand side
is well defined since $\Par(-\nu)(X,Y)$ is smooth relatively to the variable $Y\in\ggot$.
\end{prop}

\begin{proof}
The proof which is done in \cite{pep-vergneIII} follows from the
relation
\begin{equation}\label{eq:omega-P-M-nu}
\omega_P=Q^*(\omega_M)-\pi^*_2(\nu).
\end{equation}
\end{proof}

\bigskip

We now analyze the term
$$
\index^{K\times G,P}_c([Q^*\sigma])\vert_{(1,1)}(X,Y)=(2i\pi)^{-{\rm dim} P}\int_{\mathbf{T}^*P}
\Ag(P)^2\ch_c(Q^*\sigma,\omega_P)(X,Y).
$$

An easy computation gives that $\Ag(P)^2(X,Y)=j_{\kgot}(X)^{-1} q^* \Ag(M)^2(Y)$, with \break
$j_\kgot(X)=\det_{\kgot}\left(\frac{\e^{\ad(X)/2}-\e^{-\ad(X)/2}}{\ad(X)}\right)$.
If we use Proposition \ref{prop:ch-c-P-M}, we see that
\begin{eqnarray}\label{eq:calcul-indice-Q-sigma-1}
\lefteqn{\index^{K\times G,P}_c([Q^*\sigma])\vert_{(1,1)}(X,Y)  }\nonumber\\
&=&\frac{(2i\pi)^{-{\rm dim} P}}{ j_{\kgot}(X)}\int_{\mathbf{T}^*P}\!\!
\pi_1^*\!\circ \! q^*\! \left(\Ag(M)^2\ch_c(\sigma,\omega_M)\right)\!(Y)\wedge \pi_2^*\Par(-\nu)(X,Y)\nonumber\\
&=&\frac{(2i\pi)^{-{\rm dim} P}}{ j_{\kgot}(X)}\int_{\mathbf{T}^*_K P}
q^*\left(\Ag(M)^2\ch_c(\sigma,\omega_M)\right)\!(Y)\wedge \int_{\kgot^*}\Par(-\nu)(X,Y).
\end{eqnarray}

Let us compute the integral $\int_{\kgot^*}\Par(-\nu)(X,Y)$.

We choose a $K$-invariant scalar product on $\kgot$  and an
orthonormal basis $E^1,\ldots,E^r$ of $\kgot$, with dual basis
$E_1,\ldots, E_r$ : we write $X=\sum_k X_k E^k$ for $X\in\kgot$, and $\xi=\sum_k\xi_k E_k$ for
$\xi\in\kgot^*$. Let $\theta_k=\langle E_k,\theta\rangle$ be the
$1$-forms on $P$ associated to the connection one form. Let ${\rm
vol}(K,dX^o)$ be the volume of $K$ computed with the Haar measure
compatible with the volume form $dX^o= dX_1\ldots dX_r$.

We have $d\nu=\sum_k \xi_k d\theta_k + d\xi_k\theta_k$, and (\ref{eq:omega-P-M-nu})
gives that
$$
\left(d\omega_P\right)^{\dim P}= q^*\left(d\omega_M\right)^{\dim M} \wedge \theta_r\cdots\theta_1\wedge
\pi^*_2\left(d\xi_1\cdots d\xi_r\right).
$$
So, in the integral (\ref{eq:calcul-indice-Q-sigma-1}), the vector space $\kgot^*$ is oriented by the
volume form $d\xi^o= d\xi_1\cdots d\xi_r$, and $\T^*_KP$ is oriented by
$q^*\left(d\omega_M\right)^{\dim M} \wedge\theta_r\cdots\theta_1$.

Let $\Theta=d\theta+\frac{1}{2}[\theta,\theta]\in \Acal^2(P)\otimes
\kgot$ be the curvature of $\theta$.  The equivariant curvature of $\theta$ is
$$
\Theta(Y)=\mu(Y)+\Theta.
$$
Then $\Theta(Y)\in\Acal(P)\otimes \kgot$ is horizontal, and the element
$\Theta\in \Acal^2(P)\otimes\kgot$ is nilpotent. If $\varphi$ is a $\f$ function on
$\kgot$, then $\varphi(\Theta(Y))$ is computed via the Taylor series expansion at $\mu(Y)(p)$ and
$\varphi(\Theta(Y))$ is a horizontal form on $P$ which depends smoothly and $G$-equivariantly
of $Y\in\ggot$. When $\varphi\in\f(\kgot)$ is $K$-invariant, the form $\varphi(\Theta(Y))$ is basic :
hence we can look at it as a differential form on $M$ which depends smoothly and $G$-equivariantly of $Y\in\ggot$.

\begin{defi}\label{def:delta-Theta}
Let $\delta(X-\Theta(Y))$ be the $K\times G$-equivariant form on $P$ defined by the relation
$$
\int_{\kgot\times\ggot}\delta(X-\Theta(Y)) \varphi(X,Y)dXdY:={\rm vol}(K,dX)\int_{\ggot}\varphi(\Theta(Y),Y)dY,
$$
for any $\varphi\in\f(\kgot\times\ggot)$ with compact support. Here ${\rm vol}(K,dX)$ is the volume of $K$ computed
with the Haar measure compatible with $dX$.
\end{defi}

One sees that $\delta(X-\Theta(Y))$ is a $K\times G$-equivariant form on $P$ which depends smoothly
of the variable $Y\in\ggot$.

\begin{lem}\label{lem:integrale-P}Let $\kgot^*$ be oriented by the volume form $d\xi^o=d\xi_1\cdots d\xi_r$.
Then
$$
\int_{\kgot^*}\Par(-\nu)(X,Y)=(2i\pi)^{\dim K} \delta(X-\Theta(Y))
 \frac{\theta_r\cdots\theta_1}{{\rm vol}(K,dX^o)}.
$$
\end{lem}

\begin{proof}
Take $\chi_{\kgot^*}(\xi)= g(\|\xi\|^2)$ where $g\in\f_c(\Rbb)$ is
equal to $1$ in a neighborhood of $0$. Let $\varphi\in\f_c(\kgot)$
and let $\widehat{\varphi}(\xi)=\int_\kgot \e^{i\langle
 \xi,X\rangle}\varphi(X)dX^o$ be its Fourier transform relatively to $dX^o$.

To compute the integral over the fiber $\kgot^*$ of
$\Par(-\nu)(X,Y)$, only the highest exterior degree term in $d\xi$
will contribute to the integral. This term comes only from the
term $d\chi_{\kgot^*}\beta(-\nu)(X,Y)$ in
$\Par(-\nu)(X,Y):=\chi_{\kgot^*}+d\chi_{\kgot^*}\beta(-\nu)(X,Y)$.
We compute

\begin{eqnarray*}
  \int_\kgot\left(\int_{\kgot^*}\Par(-\nu)(X,Y)\right) \varphi(X)dX^o &=&
  \int_{\kgot^*}\left(\int_{\kgot}\Par(-\nu)(X,Y) \varphi(X)dX^o\right) \\
   &=& \int_{\kgot^*} d\chi_{\kgot^*}(i\nu) \left(\int_0^\infty
\e^{-it (d\nu+\langle \xi,\mu(Y)\rangle)}\widehat{\varphi}(t \xi) dt\right)\\
 &=& \int_0^\infty \underbrace{\left( \int_{\kgot^*} d\chi_{\kgot^*}(i\nu)
\e^{-it (d\nu+\langle \xi,\mu(Y)\rangle)}\widehat{\varphi}(t \xi)
\right)}_{I(t)}dt.
\end{eqnarray*}
Since $d\nu=\sum_k \xi_k d\theta_k + d\xi_k\theta_k $, the differential form $d\chi_{\kgot^*}(i\nu)\,\e^{-it d\nu}$ is
equal to
$$
2i \,
g'(\|\xi\|^2)(\sum_j \xi_j d\xi_j)(\sum_k \xi_k\theta_k)
\prod_l(1-itd\xi_l\theta_l) \e^{-it \langle \xi,d\theta \rangle},
$$
and its component $[d\chi_{\kgot^*}(i\nu)\,\e^{-it d\nu}]_{\rm max}$ of highest exterior degree  in $d\xi$ is
\begin{eqnarray*}
[d\chi_{\kgot^*}(i\nu)\,\e^{-it d\nu}]_{\rm max}
&=& -2(-i)^r t^{r-1} g'(\|\xi\|^2)\|\xi\|^2 \prod_j \left(d\xi_j\wedge\theta_j\right)
 \e^{-it \langle \xi,d\theta \rangle} \\
&=&-2(i)^r t^{r-1}\,\theta_r\cdots\theta_1\, g'(\|\xi\|^2)\|\xi\|^2  \e^{-it \langle \xi,d\theta \rangle}d\xi^o .
\end{eqnarray*}

So for $t>0$ we have
\begin{eqnarray*}
I(t)&=&-2(i)^r t^{r-1}\theta_r\cdots
 \theta_1\left(\int_{\kgot^*}g'(\|\xi\|^2)\|\xi\|^2
 \e^{-it \langle \xi,d\theta+\mu(Y) \rangle}\widehat{\varphi}(t \xi) d\xi^o\right)\\
 &=&(i)^r\theta_r\cdots\theta_1 \left(\int_{\kgot^*}\left[(-2
 g'(\hbox{$\frac{\|\xi\|^2}{t^2}$})\hbox{$\frac{\|\xi\|^2}{t^3}$}\right]
 \e^{-i \langle \xi,d\theta+\mu(Y) \rangle}\widehat{\varphi}( \xi) d\xi^o\right)\\
 &=&(i)^r\theta_r\cdots \theta_1\frac{d}{dt}\left(\int_{\kgot^*}g(\hbox{$\frac{\|\xi\|^2}{t^2}$})
 \e^{-i \langle\xi,d\theta+\mu(Y) \rangle}\widehat{\varphi}(\xi) d\xi^o\right).
\end{eqnarray*}

Finally $\int_\kgot\left(\int_{\kgot^*}\Par(-\nu)(X,Y)\right) \varphi(X)dX^o$ is equal to
\begin{eqnarray*}
\int_0^\infty I(t)dt &=&(i)^r \theta_r\cdots\theta_1\left(\int_{\kgot^*}\e^{-i \langle \xi,d\theta+\mu(Y)
\rangle}\widehat{\varphi}( \xi)d\xi^o\right)\\
&=&(2i\pi)^r \ \theta_r\cdots\theta_1\ \varphi(d\theta+\mu(Y))\\
&=&(2i\pi)^r \ \theta_r\cdots\theta_1\ \varphi(\Theta+\mu(Y)) \\
&=&(2i\pi)^r\left(\int_{\kgot}\delta(X-\Theta(Y))\varphi(X)dX^o\right)\frac{\theta_r\cdots\theta_1}{{\rm vol}(K,dX^o)}.
\end{eqnarray*}
\end{proof}

 The last lemma shows that $\index^{K\times G,P}_c([Q^*\sigma])\vert_{(1,1)}(X,Y)$
is equal to
\begin{eqnarray}\label{eq:calcul-indice-Q-sigma-2}
\lefteqn{ \frac{(2i\pi)^{-{\rm dim} M}}{ j_{\kgot}(X)}\int_{\mathbf{T}^*_K P}
q^*\left(\Ag(M)^2\ch_c(\sigma,\omega_M)\right)\!(Y)\, \delta(X-\Theta(Y))\frac{\theta_r\cdots\theta_1}{{\rm vol}(K,dX^o)}}
\nonumber\\
&=& \frac{(2i\pi)^{-{\rm dim} M}}{ j_{\kgot}(X)}\int_{\mathbf{T}^*M}
\Ag(M)^2(Y)\ch_c(\sigma,\omega_M)(Y)\,\delta_o(X-\Theta(Y)).
\end{eqnarray}

Here $\delta_o(X-\Theta(Y))$ denotes the closed $K\times G$-equivariant form $M$ defined by the
relation
$$
\int_{\kgot}\delta_o(X-\Theta(Y))\varphi(X)dX={\rm vol}(K,dX)\overline{\varphi}(\Theta(Y))
$$
for any $\varphi\in\f_c(\kgot)$. Here $\overline{\varphi}(X):= {\rm vol}(K,dk)^{-1}\int_K\varphi(kX)dk$ is the $K$-invariant function
obtained by averaging $\varphi$.

\medskip

Now we analyze the right hand side of (\ref{eq:preuve-s-s'}) at $(s,s')=(1,1)$. Here  the Chern class
$\chs(\sigma_{\tau^*})(Y)$ is equal to $\chs(\sigma)(Y)\ch(\Vcal_{\tau^*})(Y)$ where the equivariant
Chern character $\ch(\Vcal_{\tau^*})(Y)$ is represented by $\tr(\e^{\Theta(Y)},\tau^*)$. Hence
$\chc(\sigma_{\tau^*},\omega_M)(Y)=\chc(\sigma,\omega_M)(Y)\tr(\e^{\Theta(Y)},\tau^*)$.
So the generalized function \break
$\sum_{\tau\in \hat{K}} \tr(\e^X,\tau) \index^{G,M}_c([\sigma_{\tau^*}])\vert_{1}(Y)$ is equal to
\begin{equation}\label{eq:calcul-indice-Q-sigma-3}
(2i\pi)^{-{\rm dim} M} \int_{\mathbf{T}^*M}
\Ag(M)^2(Y)\ch_c(\sigma,\omega_M)(Y)\,\Xi(X,\Theta(Y))
\end{equation}
where $\Xi(X,X')$ is a generalized function on a neighborhood of $0$ in $\kgot\times\kgot$
defined by the relation  $\Xi(X,X')=\sum_{\tau\in \hat{K}} \tr(\e^X,\tau) \tr(\e^{X'},\tau^*)$.

The Schur orthogonality relation shows that
$$
\Xi(X,X')=j_{\kgot}(X)^{-1}\delta_o(X-X').
$$
In other words, $\Xi(X,X')$ is smooth relatively to $X'$ and for any $\varphi\in\f(\kgot)^K$ which is
supported in a small neighborhood of $0$, we have ${\rm vol}(K,dX) \varphi(X')$
$=\int_{\kgot}\Xi(X,X')j_{\kgot}(X)\varphi(X)dX$.

Finally, we have proved that the generalized functions (\ref{eq:calcul-indice-Q-sigma-2})
and (\ref{eq:calcul-indice-Q-sigma-3}) coincides: the proof of
(\ref{eq:preuve-s-s'}) is then completed for $(s,s')=(1,1)$.

\subsection{Multiplicative property}

We consider the setting of Subsection \ref{subsec:exterior-product}. We will check that
the cohomological index satisfies the {\em Mutiplicative property}
(see Theorem \ref{theo:multiplicative-property}).

Let $M_1$ be a compact $K_1\times K_2$-manifold, and let $M_2$ be a $K_2$-manifold. We consider the
product $M:=M_1\times M_2$ with the action of $K:=K_1\times K_2$.

\begin{theo}[Multiplicative property]
For any $[\sigma_1]\in \KK_{K_1\times K_2}(\T^*_{K_1} M_1)$ and
any $[\sigma_2]\in\KK_{K_2}(\T^*_{K_2} M_2)$ we have
\begin{equation}\label{eq:multiplicative-c-index}
\index^{K,M}_c([\sigma_1]\odot_{\rm ext}[\sigma_2])
=\index^{K_1,K_2,M_1}_c([\sigma_1])\index^{K_2,M_2}_c([\sigma_2]).
\end{equation}
The product on the right hand side  of (\ref{eq:multiplicative-c-index}) is well defined since
$\index^{K_1,K_2,M_1}_c([\sigma_1])$ is a generalized function on $K_1\times K_2$
which is smooth relatively to $K_2$ (see Lemma \ref{lem:index-smooth-H}).
\end{theo}

\begin{proof} Let $\sigma_1$ be a morphism on $\T^*M_1$, which is $K_1\times K_2$-equivariant and
$K_1$-transversally elliptic. Let $\sigma_2$ be a morphism on $\T^*M_2$, which is $K_2$-transversally elliptic.
The morphism $\sigma_2$ can be chosen so that it is almost homogeneous of degree $0$. Then the product
$\sigma:=\sigma_1\odot_{\rm ext}\sigma_2$ is $K$-transversally elliptic morphism on $\T^*M$, and
$[\sigma]=[\sigma_1]\odot_{\rm ext}[\sigma_2]$ in $\KK_{K}(\T^*_{K} M)$.

We have to show that for any $s=(s_1,s_2)\in K_1\times K_2$, we have
\begin{eqnarray}\label{eq:index-M-12}
\lefteqn{\index^{K,M}_c([\sigma])\vert_s(Y_1,Y_2)=}\\
&&\index^{K_1,K_2,M_1}_c([\sigma_1])\vert_{s_1}(Y_1,Y_2)\index^{K_2,M_2}_c([\sigma_2])\vert_{s_2}(Y_2)\nonumber
\end{eqnarray}
for $(Y_1,Y_2)$ in a neighborhood of $0$ in $\kgot_1(s_1)\times\kgot_2(s_2)$. We
conduct the proof only for $s$ equal to the identity $e$, as the
proof for $s$ general is entirely similar.

For $k=1,2$, let  $\pi_k: \T^*M\to \T^*M_k$ be the
projection. The Liouville one form $\omega$ on $\T^*(M_1\times
M_2)$ is equal to $\pi_1^*\omega_1+\pi_2^*\omega_2$, where $\omega_k$
is the Liouville one form on $\T^*M_k$.

We have three index formulas:
\begin{eqnarray*}
\index^{K,M}_c([\sigma])\vert_e(X_1,X_2)&:=&(2i\pi)^{-\dim M}
\int_{\mathbf{T}^*M}\Ag(M)^2\chc(\sigma,\omega)(X_1,X_2),\\
\index^{K,M_1}_c([\sigma_1])\vert_e(X_1,X_2)&:=&(2i\pi)^{-\dim M_1}
\int_{\mathbf{T}^*M_1}\Ag(M_1)^2\chcf(\sigma_1,\omega_1)(X_1,X_2),\\
\index^{K_2,M_2}_c([\sigma_2])\vert_e(X_2)&:=&(2i\pi)^{-\dim M_2}
\int_{\mathbf{T}^*M_2}\Ag(M_2)^2\chc(\sigma_2,\omega_2)(X_2).
\end{eqnarray*}

Following Remark \ref{rem:chc-smooth-H}, $\chcf(\sigma_1,\omega_1)(X_1,X_2)$ denotes a closed
equivariant form with compact support which represents the class  $\chc(\sigma_1,\omega_1)$, and
which is {\em smooth} relatively to $X_2\in\kgot_2$.

It is immediate to check that $\Ag(M)^2(X_1,X_2)=\Ag(M_1)^2(X_1,X_2)\Ag(M_2)^2(X_2)$. Hence
Equality (\ref{eq:index-M-12}) follows from the following identity in
$\Hcal^{-\infty}_{c}(\kgot_1\times \kgot_2,\T^*M)$ that we proved in \cite{pep-vergneIII}:
$$
\pi_1^*\chcf(\sigma_1,\omega_1)(X_1,X_2)\wedge
\pi_2^*\chc(\sigma_2,\omega_2)(X_2)=\chc(\sigma,\omega)(X_1,X_2).
$$

\end{proof}

\subsection{Normalization conditions}

\subsubsection{Atiyah symbol}\label{subsec:Atiyah-cohomology}

Let $V:= \Cbb_{[1]}$ be equipped with the canonical action of
$S^1$. The Atiyah symbol $\sigma_\at$ was introduced in Subsection
\ref{subsec:Atiyah-analytic} :  it is a $S^1$-transversally
elliptic symbol on $V$. It is the first basic example of a
``pushed'' symbol (see Subsection \ref{subsec:pushed}).

We consider on $V$ the euclidean metric $(v,w)=\Re(v\overline{w})$
: it gives at any $v\in V$ identifications $\T_v V\simeq \T^*_v
V\simeq \Cbb_{[1]}$. So in this example we will make no
distinction between vectors fields and $1$-forms on $V$.  Let
$\kappa(\xi_1)=i\xi_1$ be the vector field on $V$ associated to
the action of $S^1$ : $\kappa=-VX$ where $X=i\in {\rm Lie}(S^1)$.

Let $\sigma_V$ be the symbol on the complex vector space $V$: at any $(\xi_1,\xi_2)\in\T^* V$,
$\sigma_V(\xi_1,\xi_2): \wedge^0 V\to \wedge^1 V$ acts by multiplication by $\xi_2$. We see then that
$$\sigma_\at(\xi_1,\xi_2)=\sigma_V(\xi_1,\xi_2+\kappa(\xi_1)).$$
The symbol $\sigma_\at$ is obtained by ``pushing'' the symbol $\sigma$ by the vector field $\kappa$.

We can attached to the one form $\kappa$, the equivariant form $\Par(\kappa)$ which is defined on $V$,
and localized near $\{\kappa=0\}=\{0\}\subset V$. Since the support of $\sigma_V$ is the zero section,
the equivariant Chern character $\chs(\sigma_V)$ is an equivariant form on $\T^*V$ which is compactly
supported in the fibers of $p:\T^*V\to V$. Then the product $\chs(\sigma_V)p^*\Par(\kappa)$ defines an
equivariant form with compact support on $\T^*V$.

Here we will use the relation (see Proposition \ref{prop:chc-sigma-kappa})
\begin{equation}\label{eq:ch-sigma-kappa}
\chs(\sigma_\at)\Par(\omega)=\chs(\sigma_V)\,p^*\Par(\kappa)\quad {\rm in} \quad \Hcal^{-\infty}_c(\kgot,\T^*V).
\end{equation}

Using (\ref{eq:ch-sigma-kappa}), we now compute the
cohomological index of the Atiyah symbol.

\begin{prop}\label{indexca} We have
$$
\mathbf{[N3]}\qquad \index_c^{S^1,V}([a])(\e^{i\theta})=-\sum_{n=1}^{\infty}\e^{in\theta}.
$$
\end{prop}

\begin{proof}
We first prove  the equality above when $s=\e^{i\theta}$ is not
equal to $1$. Then, near $s$, the generalized function
$-\sum_{n=1}^{\infty}\e^{in\theta}$ is analytic  and given by
$-\frac{s}{1-s}$.

Now, at a point $s\in S^1$ different from $1$, the fixed point set
$V(s)$ is $\{0\}$. The character $\ch_s(\Ecal)$ is $(1-s)$, and
the form $D_s(\Ncal)$ is $(1-s)(1-s^{-1})$. Thus
$$\index_c^{S^1,V}(s)=\frac{(1-s)}{(1-s)(1-s^{-1})}=-\frac{s}{1-s}.$$
This shows the equality of both members in Proposition
\ref{indexca} on the open set $s\neq 1$ of $S^1$.

We now compute near $s=1$. Thanks to Formula (\ref{eq:ch-sigma-kappa}) we have
$$
\index_c^{S^1,V}(\sigma_\at)\vert_1(\theta)=
(2i\pi)^{-2}\int_{\mathbf{T}^*V}\Ag(V)^2(\theta)\,\chs(\sigma_V)(\theta)\, p^*\Par(\kappa)(\theta).
 $$

The Chern character with support $\chg(\sigma_V)(\theta)$ is proportional
to the $S^1$-equivariant  Thom form  of the real vector bundle
$\T^*V\to V$. More precisely, calculation already done in
\cite{pep-vergneII} shows that
$$\chs(\sigma_V)(\theta)=
(2i\pi)\frac{\e^{i\theta}-1}{i\theta}{\rm Thom}(\T^*V)(\theta).
$$
However the symplectic orientation on $\T^*V\simeq\Cbb^2$ is the opposite of the
orientation given by its complex structure. Now
$$\Ag(V)^2(\theta)=\frac{(i\theta)(-i\theta)}{(1-\e^{i\theta})(1-\e^{-i\theta})}.$$
Thus we obtain
$$
\index_c^{S^1,V}(\sigma_\at)\vert_1(\theta)=\frac{-i\theta}{(1-\e^{-i\theta})}
\frac{1}{2i\pi}\int_V \Par(\kappa)(\theta).
$$
As $\frac{(1-\e^{-i\theta})}{-i\theta}=-\int_{-1}^0\e^{ix\theta}dx$,
we see that
$\frac{(1-\e^{-i\theta})}{-i\theta}(-\sum_{n=1}^{\infty}\e^{in\theta})
=\int_{0}^{\infty}\e^{ix\theta}dx$. It remains to show
\begin{equation}\label{parkappa}
\frac{1}{2i\pi}\int_V\Par(\kappa)(\theta)=\int_0^{\infty}\e^{ir\theta}dr.
\end{equation}
We have $D\kappa(\theta)=\theta(x^2+y^2)+2dx\wedge dy.$ Take a
function $g$ on $\Rbb$ with compact support and equal to $1$ on a
neighborhood of $0$. Let $\chi=g(x^2+y^2)$. Then
\begin{eqnarray*}
\Par(\kappa)(\theta)&=&\chi-id\chi\wedge
\kappa\int_0^{\infty}\e^{itD\kappa(\theta)}dt\\
&=&g(x^2+y^2)-2ig'(x^2+y^2) dx\wedge dy\int_0^{\infty}(x^2+y^2)\e^{i\theta t(x^2+y^2)}dt\\
&=&g(x^2+y^2)-2ig'(x^2+y^2)dx\wedge dy\int_0^{\infty}\e^{i\theta t}dt.
\end{eqnarray*}

Finally we obtain (\ref{parkappa}) since
$\int_V-2ig'(x^2+y^2) dx\wedge dy=2i\pi$. This completes the proof.

\end{proof}

\subsubsection{Bott symbols}

We will check here that the cohomological index satisfies the condition
$\mathbf{[N2]}$: $\qquad \index^{O(V),V}_c(\bott(V_\Cbb))=1$, for any
euclidean vector space $V$.

We have explain in Remark \ref{rem:reduction-bott-symbol} that it sufficient to prove
$\mathbf{[N2]}$ for the cases:

$\bullet$ $V=\Rbb$ with the action of the group $O(V)=\Zbb/2\Zbb$,

$\bullet$ $V=\Rbb^2$ with the action of the group $SO(V)=S^1$.

Let $V=\Rbb$ with the multiplicative action of $\Zbb_2$. We have to check
that $\index_c^{\Rbb, \Zbb_2}(\bott(\Cbb))(\epsilon)=1$ for $\epsilon\in\Zbb_2$.
When $\epsilon=1$, we have
$$
\index_c^{\Rbb, \Zbb_2}(\bott(\Cbb))(1)=(2i\pi)^{-1}\int_{\mathbf{T}^*\Rbb}
\Ag(\Rbb)^2\chc(\bott(\Cbb)).
$$
Here $\Ag(\Rbb)^2=1$. We have proved in \cite{pep-vergneII} that the class
$\chc(\bott(\Cbb))\in\Hcal_c^2(\T^*\Rbb)$ is equal to $2i\pi$ times the
Thom form of the oriented vector space of $\Rbb^2\simeq \T^*\Rbb$. Hence
$\index_c^{\Rbb, \Zbb_2}(\bott(\Cbb))(1)=1$.
When $\epsilon=-1$, the space $\T^*\Rbb(\epsilon)$ is reduced to a point.
We see that $\chc(\bott(\Cbb),\epsilon)=2$, $D_\epsilon(\Ncal)=\det(1-\epsilon)=2$.
Then $\index_c^{\Rbb, \Zbb_2}(\bott(\Cbb))(1)=\frac{\chc(\bott(\Cbb),\epsilon)}{D_\epsilon(\Ncal)}=1$.

Let $V=\Rbb^2$ with the rotation action of $S^1$. Like before
$\index_c^{\Rbb^2, S^1}(\bott(\Cbb^2))(1)$ is equal to $1$ since the Chern class $\chc(\bott(\Cbb^2))$
is equal to $(2i\pi)^2$ times the
Thom form of the oriented vector space of $\Rbb^4\simeq \T^*\Rbb^2$.
When $\e^{i\theta}\neq 1$, the space $\T^*\Rbb^2(\e^{i\theta})$ is reduced to a point.
We see that $\chc(\bott(\Cbb),\e^{i\theta})=D_{\e^{i\theta}}(\Ncal)=2(1-\cos(\theta))$.
Then $\index_c^{\Rbb^2, S^1}(\bott(\Cbb^2))(\e^{i\theta})=1$.

\section{Examples}

\subsection{Pushed symbols}\label{subsec:pushed}

Let $M$ be a $K$-manifold and $N=\T^*M$. Let $\Ecal^{\pm}\to M$ be two
$K$-equivariant complex vector bundles on $M$ and $\sigma:
p^*\Ecal^+\to p^*\Ecal^-$ be a $K$-equivariant symbol  which is
supposed to be \textbf{invertible exactly outside the zero section} : the set $\supp(\sigma)$  coincides with
the zero section of $\T^*M$.

If $M$ is compact, $\sigma$ defines an elliptic
symbol on $\T^* M$, thus a fortiori a transversally elliptic
symbol.

Here we assume $M$ non compact. Following Atiyah's strategy
\cite{Atiyah.74}, we can ``push" the symbol $\sigma$ outside the
zero section, if we dispose of a $K$-invariant  real one-form
$\kappa$ on $M$.  This construction provides  radically new
transversally elliptic symbols. We recall some definitions of
\cite{pep-vergneIII}:

\begin{defi}
Let $\kappa$ be a real $K$-invariant one-form on $M$. Define
$f_\kappa:M\to \kgot^*$ by $\langle f_\kappa(x),X\rangle=\langle \kappa(x), V_xX\rangle$.
We define the subset $C_\kappa$ of $M$ by
$C_\kappa=f_\kappa^{-1}(0)$. We call $C_\kappa$ the critical set
of $\kappa$.

\end{defi}

We define the symbol $\sigma(\kappa)$ on $M$ by
$$
\sigma(\kappa)(x,\xi)=\sigma(x,\xi+\kappa(x)),\quad {\rm for}\quad (x,\xi)\in \T^*M.
$$
Thus $\sigma(\kappa)$ is not invertible at
$(x,\xi)$ if and only if $\xi=-\kappa(x)$, and then
$(x,\xi)\in\supp(\sigma(\kappa))\cap \T^*_KM$ if $\xi=-\kappa(x)$
and $\langle\kappa(x),V_xX\rangle=0$ for all $X\in\kgot$. Thus
$$
\supp(\sigma(\kappa))\cap \T^*_KM=\{(x,-\kappa(x))\ |\  x\in
C_\kappa\}.
$$

If  $C_\kappa$ is compact, then the morphism $\sigma(\kappa)$ is
transversally elliptic.

\bigskip

Using a $K$-invariant metric on $\T M$, we can associate to a
$K$-invariant vector field $\Kcal$ on $M$ a $K$-invariant real
one-form.

\begin{exam}\label{pushS}
Let $S\in \kgot$ be a central element of $\kgot$ such that the set
of zeroes of $VS$ is compact. Then the associated  form
$\kappa_S(\bullet)=\langle VS,\bullet\rangle $ is  a $K$-invariant
real one-form such that $C_{\kappa_S}$ is compact. Indeed the
value of $f_{\kappa_S}$ on $S$ is $\|VS\|^2$, so that the set
$C_{\kappa_S}$ coincides with the fixed point set $M(S)$.
\end{exam}

\begin{defi}\label{pushsymbol}
If $\kappa$ is a $K$-invariant real one-form on $M$ such that
$C_\kappa$ is compact, the transversally elliptic symbol
$$\sigma(\kappa)(x,\xi)=\sigma(x,\xi+\kappa(x))$$ is called the
pushed symbol of $\sigma$ by $\kappa$.
\end{defi}

\begin{exam}\label{ex-pushed-Atiyah}
The Atiyah symbol is a pushed symbol defined on $M=\Rbb^2$ (see Subsection \ref{subsec:Atiyah-cohomology}).
\end{exam}

We construct as in (\ref{eq:beta-omega}) the $K$-equivariant differential
form
$$
\beta(\kappa)(X)=-i\kappa\wedge
\int_{0}^{\infty}\e^{itD\kappa(X)}dt
$$
which is defined on $M\setminus C_\kappa$. We choose a compactly supported function $\chi_{\kappa}$ on $M$
identically $1$ near $C_\kappa$. Then the $K$-equivariant form
$$
\Par(\kappa)(X)=\chi_{\kappa}+d\chi_{\kappa}\beta(\kappa)(X)
$$
defined a class in $\Hcal_{c}^{-\infty}(\kgot,M)$.

The $K$-equivariant form  $\Par(\kappa)$ is congruent to $1$ in
cohomology without support conditions. Indeed one verify
that $\Par(\kappa)=1+D((\chi_\kappa-1)\beta(\kappa))$.

Let $p:\T^*M\to M$ be the projection. We can multiply the $K$-equivariant form $p^*\Par(\kappa)(X)$ with
$C^{-\infty}$-coefficients by the $K$-equivariant form
$\chs(\sigma)(X)$. In this way, we obtain a $K$-equivariant form
with compact support on $\T^*M$.

\begin{prop}\label{prop:chc-sigma-kappa}
The $K$-equivariant form  $\chs(\sigma)p^*\Par(\kappa)$ represents the class $\chc(\sigma(\kappa),\omega)$  in
$\Hcal^{-\infty}_{c}(\kgot,\T^*M)$.
\end{prop}

\begin{proof}
By definition the class $\chc(\sigma(\kappa),\omega)$ is represented by the product $\chs(\sigma(\kappa))\Par(\omega)$
We first prove $\chs(\sigma(\kappa))\Par(\omega)=\chg(\sigma(\kappa))\Par(p^*\kappa)$ in
$\Hcal^{-\infty}_{c}(\kgot,\T^*M)$.

Indeed if $(x,\xi)\in \supp(\sigma(\kappa))$, then $\xi=-\kappa(x)$.
Thus $\langle \omega(x,\xi),v\rangle =- \langle \xi,p_*v\rangle$
$= \langle \kappa(x),p_*v\rangle $ where $v$ is any tangent vector at $(x,\xi)\in \T^*M$.
So the  $1$-forms $\omega$ and $p^*\kappa$ coincides on the support of $\sigma(\kappa)$. Thus $\chs(\sigma(\kappa))
\Par(\omega)=\chg(\sigma(\kappa))p^*\Par(\kappa)$ as consequence of (\cite{pep-vergneIII}, Corollary 3.12).

Let us prove now that $\chs(\sigma(\kappa))p^*\Par(\kappa)=\chs(\sigma)p^*\Par(\kappa)$.
Consider the family of symbols on $M$ defined by $\sigma_t(x,\xi)=\sigma(x,\xi+t\kappa(x))$
for  $t\in[0,1]$ : we have $\sigma_0=\sigma$ and $\sigma_1=\sigma(\kappa)$.

On a compact neighborhood $\Ucal$ of $C_\kappa$, the support of
$\sigma_t$ stays in the compact set $\{(x,\xi): x\in \Ucal, \xi=-t\kappa(x)\}$
when $t$ varies between $0$ and $1$. It follows
from (\cite{pep-vergneIII}, Theorem 3.11)  that all the classes
$\chs(\sigma_t)p^*\Par(\kappa),\ t\in[0,1]$ coincides in $\Hcal_{c}^{-\infty}(\kgot,\T^*M)$.
\end{proof}

\medskip

Similarly for any $s\in K$, we consider the restriction $\kappa_s$
of the form $\kappa$ to $M(s)$. We finally obtain the following
formula:

\begin{theo}\label{indexpushed}
For any $s\in K$ and $X\in \kgot(s)$ small, the cohomological
index $\index_c^{K,M}(\sigma(\kappa))\vert_s(Y)$ is given on
$\kgot(s)$ by the integral formula:
$$
\int_{\mathbf{T}^*M(s)}\Lambda_s(Y) \,\Par(\kappa_s)(Y) \, \chs(\sigma,s)(Y).
$$
In particular, when $s=1$ we get
$$
\index_c^{K,M}(\sigma(\kappa))\vert_e(X)=(2i\pi)^{-\dim M}\int_{\mathbf{T}^*M}\Ag(M)^2(X) \,
\Par(\kappa)(X) \, \chs(\sigma)(X).
$$
\end{theo}

\medskip

An interesting situation is when the manifold $M$ is oriented, and is equipped with a $K$-invariant Spin structure.
Let $\Scal_M\to M$ be the corresponding spinor bundle. We consider the $K$-invariant symbol
$\sigma_{spin}:p^*\Scal_M^+\to p^*\Scal_M^-$. It support is exactly the zero section of the cotangent bundle.
For any invariant $1$-form $\kappa$ on $M$ such that $C_\kappa$ is compact we consider the
transversally elliptic symbol $\sigma_{\spin}(\kappa)$.

We have proved in \cite{pep-vergneII} that
$$
\chs(\sigma_{spin})(X)=(2i\pi)^{\dim M} \hat{A}(M)^{-1}(X)\, {\rm Thom}(\T^*M)(X).
$$
Hence Theorem \ref{indexpushed} tell us that
$$
\index_c^{K,M}(\sigma_{spin}(\kappa))\vert_e(X)= \int_{M}\Ag(M)(X)\,\Par(\kappa)(X).
$$

\subsection{Contact manifolds}

The following geometric example is taken from \cite{Fitzpatrick07}.

Let $M$ be a compact manifold of dimension $2n+1$. Suppose that $M$ carries a {\em contact} $1$-form $\alpha$ ;
that is, $E=\ker(\alpha)$ is a hyperplane distribution of $\T M$, and the restriction of the $2$-form
$d\alpha$ to $E$ is symplectic. The existence of a Reeb vector field $\mathbf{Y}$ on $M$ gives canonical
decompositions $\T M=E\oplus \Rbb \mathbf{Y}$ and $\T^*M=E^*\oplus E^0$ with $E^0=\Rbb\alpha$.

Let $J$ be a $K$-invariant complex structure on the bundle $E$ which is compatible with the
symplectic structure $d\alpha$. We equipped the bundle $E^*$ with the complex structure $J^*$
defined by $J^*(\xi):=\xi\circ J$ for any cotangent vector $\xi$.
We note that the complex bundle $(E^*,J^*)$ is the complex dual of the vector bundle $(E,J)$.

We consider the $\Zbb_2$-graded complex vector bundle $\Ecal:=\wedge_{J^*} E^*$.
The Clifford action defines a bundle map
$\clif: E^*\to \End_\Cbb(\Ecal)$. We consider now the symbol on $M$
$$
\sigma_b: p^*(\Ecal^+)\to p^*(\Ecal^-)
$$
defined by $\sigma_b(x,\xi)=\clif(\xi')$ where $\xi'$ is the projection of $\xi\in \T^*M$ on $E^*$.

We see that the support of $\sigma_b$ is equal to $E^0\subset \T^*M$ : $\sigma_b$ is not an elliptic symbol.

Let $K$ be a compact Lie group acting on $M$, which leaves $\alpha$ invariant. Then $E,E^*$ are
$K$-equivariant complex vector bundles, and the complex struture $J$ can be chosen $K$-invariant.
The morphism $\sigma_b$ is then $K$-equivariant.

We suppose for the rest of this section that
\begin{equation}\label{hypothese-contact}
E^0\cap \T^*_K M=\ {\rm zero\ section\ of}\ \ \T^*M
\end{equation}
It means that for any $x\in M$, the map $f_\alpha(x):X\mapsto \alpha_x(V_xX)$ is not the zero map.
Under this hypothesis the symbol $\sigma_b$ is transversally elliptic.

Under the hypothesis (\ref{hypothese-contact}), we can define the following closed equivariant form on
$M$ with $\fgene$-coefficients
$$
\Jcal_\alpha(X):=\alpha\int_\Rbb \e^{itD\alpha(X)}dt.
$$

For any $\varphi\in\f_c(\kgot)$, the expression $\int_\kgot\Jcal_\alpha(X)\varphi(X)dX:=
\alpha\int_\Rbb \e^{itd\alpha}\widehat{\varphi}(t f_\alpha)dt$ is a well defined differential form on $M$
since the map $f_\alpha:M\to\kgot^*$ as an empty $0$-level set.

Let ${\rm Todd}(E)(X)$ be the equivariant Todd class of the complex vector bundle $(E,J)$. We have the following

\begin{theo}[\cite{Fitzpatrick07}]
For any $X\in\kgot$ sufficiently small,
$$
\index^{K,M}_c(\sigma_b)(\e^X)=(2i\pi)^{-n}\int_M {\rm Todd}(E)(X)\,\Jcal_\alpha(X).
$$
\end{theo}

\begin{proof}
Consider the equivariant form  with compact support $\chs(\sigma_b)\Par(\omega)$. The Chern form
$\chs(\sigma_b)$ attached to the complex vector bundle $E^*$ is computed in \cite{pep-vergneIII}
as follows. Let ${\rm Thom}(E^*)(X)$ be the equivariant Thom
form, and let ${\rm Todd}(E^*)(X)$ be the equivariant Todd form. We have proved in \cite{pep-vergneIII}, that
\begin{equation}\label{eq:ch-sigma-contact}
\chs(\sigma_b)=(2i\pi)^{n}\,{\rm Todd}(E^*)(X)^{-1}{\rm Thom}(E^*)(X).
\end{equation}

Let $[\Rbb]$ be the trivial vector bundle over $M$. We work trought the isomorphism $E^*\oplus[\Rbb]\simeq \T^*M$
who sends $(x,\xi,t)$ to $(x,\xi+t\alpha(x))$. We consider the invariant $1$-form $\lambda$ on
$E^*\oplus[\Rbb]\simeq \T^*M$ defined by
$$
\lambda= -\underline{t} \ p^*(\alpha)
$$
Here $p:E^*\oplus[\Rbb]\to M$ is the projection, and $\underline{t}$ denotes the function that sends $(x,\xi,t)$ to $t$.

It is easy to check that the form $\lambda$ and the Liouville one form $\omega$ are equal on the support of $\sigma_b$.
Thus
$$
\chs(\sigma_b)\Par(\omega)=\chs(\sigma_b)\Par(\lambda), \quad {\rm in}\quad \Hcal^{-\infty}_c(\kgot,\T^*M),
$$
as consequence of (\cite{pep-vergneIII}, Corollary 3.12). We have then
$$
\index^{K,M}_c(\sigma_b)(\e^X)=(2i\pi)^{-\dim M}\int_{\mathbf{T}^*M} \Ag(M)^2(X)\chs(\sigma_b)(X)\Par(\lambda)(X).
$$
The integral of $\chs(\sigma_b)(X)\Par(\lambda)(X)$ on the fibers of $\T^*M$ is then equal to the
product
$$
\left(\int_{E^*\, fiber}\chs(\sigma_b)(X)\right)\left(\int_{\Rbb}\Par(\lambda)(X)\right)
$$
If we uses (\ref{eq:ch-sigma-contact}), we see that the  integral $\int_{E^*\, fiber}\chs(\sigma_b)(X)$  is equal to
$(2i\pi)^{n}{\rm Todd}(E^*)(X)^{-1}$. A small computation gives that $\int_{\Rbb}\Par(\lambda)(X)$
is equal to $(2i\pi)\Jcal_\alpha(X)$. The proof is now completed since
$\Ag(M)^2(X){\rm Todd}(E^*)(X)^{-1}={\rm Todd}(E)(X)$.
\end{proof}

\end{document}